\documentclass[pdflatex,sn-mathphys-num]{sn-jnl}

\usepackage{amsmath,amssymb,amsthm}
\usepackage{algorithm}
\usepackage{algpseudocode}
\usepackage{multirow} 
\usepackage{lmodern}
\usepackage{float}

\newtheorem{theorem}{Theorem}[section]
\newtheorem{proposition}[theorem]{Proposition}

\newtheorem{algo}[theorem]{Algorithm}

\theoremstyle{definition}

\newtheorem{remark}[theorem]{Remark}
\newtheorem{definition}[theorem]{Definition}
\newtheorem{example}[theorem]{Example}

\allowdisplaybreaks
\def\QQ{\mathbb{Q}}
\def\AA{\mathbb{A}}
\def\FF{\mathbb{F}}
\def\BB{\mathbb{B}}
\def\II{\mathbb{I}}
\def\NN{\mathbb{N}}
\def\ZZ{\mathbb{Z}}
\def\TT{\mathbb{T}}

\newcommand{\M}{{\mathfrak{M}}}

\let\epsilon=\varepsilon

\def\phi{{\varphi}}
\let\Psi=\varPsi
\let\Phi=\varPhi
\let\theta=\vartheta
\let\rho=\varrho

\def\LT{\mathop{\rm LT}\nolimits}

\def\NR{\mathop{\rm NR}\nolimits}

\def\Mat{\mathop{\rm Mat}\nolimits}

\def\Supp{\mathop{\rm Supp}\nolimits}

\def\Ker{\mathop{\rm Ker}\nolimits}
\def\rk{\mathop{\rm rk}\nolimits}

\newcommand{\Lin}{\mathop{\rm Lin}\nolimits}
\newcommand{\edim}{\mathop{\rm edim}\nolimits}

\newcommand{\cms}{\texttt{Cryp\-to\-Mi\-ni\-Sat}}
\newcommand{\kissat}{\texttt{Kis\-sat}}
\newcommand{\bosphorus}{\texttt{Bos\-pho\-rus}}

\def\Gstandard{{G_{\rm std}}}
\def\G6{{G_{1-6}}}

\def\k{{\tt k}}
\def\p{{\tt p}}
\def\c{{\tt c}}
\def\rk{{\tt rk}}
\def\si{{\tt si}}
\def\so{{\tt so}}
\def\rks{{\tt rks}}

\def\tsum_#1^#2{{\textstyle\sum\limits_{#1}^{#2}}}

\def\cocoa{\mbox{\rm
  C\kern-.13em o\kern-.07 em C\kern-.13em o\kern-.15em A}}
\def\apcocoa{\mbox{\rm
A\kern-0.13em p\kern -0.07em C\kern-.13em o\kern-.07 em C\kern-.13em
o\kern-.15em A}}

\begin{document}

\title[Efficiently Checking Separating Indeterminates]{Efficiently Checking Separating Indeterminates}

\author[1]{\fnm{Bernhard} \sur{Andraschko}}\email{Bernhard.Andraschko@uni-passau.de}
\author*[1]{\fnm{Martin} \sur{Kreuzer}}\email{Martin.Kreuzer@uni-passau.de}
\author[2]{\fnm{Le Ngoc} \sur{Long}}\email{lelong@hueuni.edu.vn}

\affil*[1]{\orgdiv{Fakult\"at f\"ur Informatik und Mathematik}, \orgname{Universit\"at Passau}, 
   \orgaddress{\street{Innstra{\ss}e~33}, \postcode{D-94032} \city{Passau}, \country{Germany}}}
\affil*[2]{\orgdiv{University of Education}, \orgname{Hue University}, 
   \orgaddress{\street{34 Le Loi Street}, \city{Hue City}, \country{Vietnam}}}

\abstract{In this paper we continue the development of a new technique for computing
elimination ideals by substitution which has been called $Z$-separating re-embeddings.
Given an ideal $I$ in the polynomial ring $K[x_1,\dots,x_n]$ over a field $K$,
this method searches for tuples $Z=(z_1,\dots,z_s)$ of indeterminates with the property
that $I$  contains polynomials of the form $f_i = z_i - h_i$ for
$i=1,\dots,s$ such that no term in $h_i$ is divisible by an indeterminate in $Z$.
As there are frequently many candidate tuples $Z$,
the task addressed by this paper is to efficiently check whether
a given tuple $Z$ has this property. We construct fast algorithms
which check whether the vector space spanned by the generators of $I$ or a somewhat
enlarged vector space contain the desired polynomials $f_i$. We also extend these algorithms 
to Boolean polynomials and apply them to cryptoanalyse round reduced versions of the AES cryptosystem 
faster.}



\maketitle

%
%

\section{Introduction}
\label{sec1-intro}

Computing elimination ideals has long been one of the key tasks of computer algebra.
In algebraic geometry, they define projections of varieties or schemes,
in linear algebra and number theory, they are applied to find minimal polynomials,
and in numerous applications, they are used to help solve polynomial systems of
equations. Traditionally, elimination ideals have been computed using resultants
or using Gr\"obner bases with respect to elimination term orderings. As these techniques
become impractical when the number of indeterminates is large, a new technique,
called $Z$-separating re-embeddings, and based on trying to perform elimination  
by substitution, has been developed in~\cite{KLR1,KLR2,KLR3,KR3}. This technique proceeds
as follows.

Let $P=K[x_1,\dots,x_n]$ be a polynomial ring over a field~$K$, and let $I$ 
be an ideal in~$P$. Given a tuple of distinct indeterminates $Z=(z_1,\dots,z_s)$
in $X=(x_1,\dots,x_n)$, we say that a tuple of polynomials $(f_1,\dots,f_s)$ in~$I$ is
{\it $Z$-separating} if there exists a term ordering~$\sigma$ such that 
$\LT_\sigma(f_i) = z_i$ for $i=1,\dots,s$. If we find such a tuple of polynomials,
we can calculate a tuple of {\it coherently $Z$-separating polynomials} and 
then the elimination ideal $I\cap K[X\setminus Z]$ by potentially less costly
interreductions and substitutions. Notice that in contrast to~\cite{KLR1,KLR2,KLR3}
we do not assume $I \subseteq \langle x_1,\dots,x_n\rangle$ here, as this hypothesis 
is not needed for our algorithms to work.

Fast methods for determining good candidate tuples~$Z$ for which a $Z$-separating
tuple of polynomials might exist in~$I$ were developed in~\cite{KLR2,KLR3,KR3}. Here we 
study the second phase of the method, namely to determine quickly whether the ideal~$I$
contains a $Z$-separating tuple of polynomials for a given~$Z$. To answer this
question in full generality, one would have to compute a Gr\"obner basis of~$I$, the very task
that we deem infeasible and want to avoid. Hence we will content ourselves with an algorithm which
checks very quickly if the vector space spanned by the given generators of~$I$,
or a somewhat larger, easily calculated vector space, contains a $Z$-separating tuple.
These algorithms are fast enough to allow us to scan many hundred candidate tuples~$Z$
in a matter of seconds, and are therefore suitable for applications involving larger
polynomial systems. For instance, in the final section we apply them to speed up 
algebraic attacks on the AES-128 cryptosystem.

Now let us describe the contents of the paper in more detail. In Section~\ref{sec2-basics}
we recall the basic definitions of separating and coherently separating tuples of polynomials,
and of separating re-embeddings. Then, in Section~\ref{sec3-check}, we introduce the main
algorithm of this paper (see Algorithm~\ref{alg-checkZ}). It allows us to discover 
very efficiently when a given polynomial ideal contains a $Z$-separating tuple.
However, it is allowed to fail in certain complicated cases. The idea of this algorithm
is to construct a weight vector $W\in \mathbb{N}^n$ such that, for every term ordering~$\sigma$
which is compatible with the grading on~$P$ given by~$W$, there exists a $Z$-separating
tuple of polynomials $(f_1,\dots,f_s)$ in~$I$ with $\LT_\sigma(f_i)=z_i$ for $i=1,\dots,s$.
Moreover, in Proposition~\ref{prop-checkZ-success}, we prove that this algorithm
succeeds exactly when the vector space spanned by the initially given system of generators of~$I$
contains a $Z$-separating tuple.

Section~\ref{sec4-checkopt} introduces an optimization of our main algorithm.
Algorithm~\ref{alg-OptimizedCheckZ} allows us to find $Z$-separating tuples of
polynomials in~$I$ in more cases than Algorithm~\ref{alg-checkZ} at the cost of
a small reduction in efficiency. The underlying idea is derived from the
border basis algorithm (see~\cite{KK}). We do not search for a
$Z$-separating tuple in the vector space~$V$ spanned by the original generators
of~$I$, but in the larger vector space $V^+ = V + x_1 V + \cdots + x_n V$.
The loss in efficiency is mostly due to the necessity to perform some linear
algebra operations (see Remark~\ref{rem-computingBasis-usingDivAlg}).

Now assume that one of these algorithms confirms the existence of a $Z$-separating
tuple in~$I$. It still remains to find an actual such tuple. This task is tackled
in Section~\ref{sec5-compute-polys}. In case the existence has been shown by Algorithm~\ref{alg-checkZ},
a straightforward base change mechanism allows us to exhibit an actual $Z$-separating
tuple of poylnomials in~$I$ (see Algorithm~\ref{alg-findSepTuple}). If the existence
has been shown using Algorithm~\ref{alg-OptimizedCheckZ}, the situation is a bit trickier:
in this case we have to mimick some linear algebra operations performed by the algorithm
on simplified polynomials using the original generators instead (see Remark~\ref{rem-findSepTupleOptim}).
We also explain how to pass from a $Z$-separating tuple to a coherently $Z$-separating tuple
(see Remark~\ref{rem-SepToCohSep}).

If we want to apply these algorithms in situations involving Boolean polynomials,
it suffices to adjust them slightly. In Section~\ref{sec6-boolpolys} we first
provide the necessary background knowledge about canonical representatives, Boolean
Gr\"obner bases, $Z$-separating and coherently $Z$-separating tuples of Boolean polynomials,
as well as Boolean $Z$-separating re-embeddings. Based on the result that the
canonical representatives of a tuple of Boolean polynomials are $Z$-separating if and only
if the tuple itself is $Z$-separating (see Proposition~\ref{prop-boolZsepConnection}), we
show how to find $Z$-separating tuples of Boolean polynomials (see Proposition~\ref{prop-checkZ-bool}).

Finally, in Section~\ref{sec7-appl}, we apply the algorithms in order to improve 
algebraic attacks on the cryptosystem AES-128. These improvements are a consequence
of better algebraic representations of the S-boxes of this cipher (see Example~\ref{ex-sbox-polys}).
They translate to more compact representations of $r$-round AES-128 in several
logic normal forms such as CNF, CNF-XOR or XNF. For one round of AES-128, these
better representations translate to meaningful speed-ups of state-of-the-art Boolean solvers
(see Example~\ref{ex-1-AES}), and for two rounds of AES, we get
faster solutions in the case when we know some key bits (see Example~\ref{ex-2-AES}).

All algorithms are illustrated by non-trival explicit examples. These examples were calculated
using implementations of the algorithms in the computer algebra systems \apcocoa 
(see~\cite{ApCoCoA}), \cocoa (see~\cite{CoCoA}), and {\tt SageMath} (see~\cite{sage}).
The source code is available upon request from the authors.
The general notation and definitions in this paper follow~\cite{KR1,KR2}.

\bigbreak
%
%

\section{Separating Tuples and Separating Re-Embeddings}
\label{sec2-basics}

In the following we let $K$ be an arbitrary field, let $P=K[x_1,\dots,x_n]$,
let~$I$ be a proper ideal in~$P$, and let $R=P/I$. Rings of this form are also called
{\it affine $K$-algebras}, since they are the affine coordinate rings of
closed subschemes of~$\mathbb{A}^n_K$. The well-known task of re-embedding those
subschemes into lower dimensional affine spaces can be phrased in terms of their coordinate
rings as follows (see~\cite{KLR1}, \cite{KLR2}).

\begin{definition}\label{def-reembed}
Let $R=P/I$ be an affine $K$-algebra.
\begin{enumerate}
\item[(a)] A $K$-algebra isomorphism $\Phi:\; R\longrightarrow P'/I'$,
where~$P'$ is a polynomial ring over~$K$ and~$I'$ is an ideal in~$P'$,
is called a \textbf{re-embedding} of $I$.

\item[(b)] A re-embedding $\Phi:\; R\longrightarrow P'/I'$ of~$I$ is called
\textbf{optimal}, if every $K$-al\-ge\-bra isomorphism $R\longrightarrow P''/I''$
with a polynomial ring $P''$ over $K$ and an ideal~$I''$ in~$P''$
satisfies $\dim(P'')\ge \dim(P')$.
\end{enumerate}
\end{definition}

In a series of papers (\cite{KLR1}, \cite{KLR2}, \cite{KLR3}) methods were developed
for finding good re-embeddings using the following technique.
Let $X =(x_1,\dots,x_n)$, and let $Z =(z_1,\dots,z_s)$ be a tuple of
pairwise distinct indeterminates in~$X$.

\begin{definition}\label{def-Zsep}
Let $I$ be a proper ideal in~$P$.
\begin{enumerate}
\item[(a)] The tuple~$Z$ is called a {\bf separating tuple of indeterminates} for~$I$,
if there exists a term ordering~$\sigma$ and a tuple of polynomials
$F=(f_1,\dots,f_s)$ with $f_i\in I$ such that $\LT_\sigma(f_i)=z_i$
for $i=1,\dots,s$.

\item[(b)] In this case we say that~$F$ is a {\bf $Z$-separating tuple} in~$I$.

\item[(c)] If, additionally, we have $f_i = z_i-h_i$ with $h_i \in K[X\setminus Z]$ for
$i=1,\dots,s$, then $F$ is called a {\bf coherently $Z$-separating tuple} in~$I$.
\end{enumerate}
\end{definition}

Notice that, in contrast to~\cite{KLR1}, \cite{KLR2}, and~\cite{KLR3}, we do not require
$I\subseteq \M = \langle x_1,\dots,x_n \rangle$ here.
Also note that the condition $\LT_\sigma(f_i)=z_i$ in this definition
implies that no other term in $\Supp(f_i)$ is divisible by~$z_i$.
When we consider the subideal $J=\langle f_1,\dots,f_s\rangle$ of~$I$,
the tuple~$F$ is a minimal monic $\sigma$-Gr\"obner basis for~$J$.
After interreducing the polynomials in~$F$,
we obtain the reduced $\sigma$-Gr\"obner basis $\widetilde{F} = (\tilde{f}_1,\dots,
\tilde{f}_s)$ of~$J$, which is then a coherently $Z$-separating tuple.

Given a coherently $Z$-separating tuple in~$I$, we can define a re-embedding of~$I$
as in the following proposition (cf.~\cite{KLR1}, \cite{KLR2}).
Here we denote by $\Lin_{\M}(I)$ the ideal
$\langle \Lin(f) \mid f \in I \rangle$, where $\Lin(f)$ denotes 
the {\it linear part} of a polynomial $f\in I$, i.e., its homogeneous component of degree~1.
Note that if~$I$ contains a polynomial with a non-zero constant term, then
$\Lin_{\M}(I) = \langle x_1,\dots,x_n\rangle_K$.
If~$I$ is contained in~$\M$ then the ideal $\Lin_{\M}(I)$ is generated
by the linear parts of a system of generators of~$I$ (see~\cite{KLR1}, Prop.~1.9).

\begin{proposition}{\bf ($Z$-Separating Re-Embeddings)}\label{prop-reembed}\\
Let~$Z$ be a tuple of distinct indeterminate in~$X$, let $Y=X\setminus Z$,
let $I$ be a proper ideal in~$P$, and let $\widetilde{F} = (\tilde{f}_1, \dots,
\tilde{f}_s)$ be a coherently $Z$-separating tuple of polynomials in~$I$.
For $i=1,\dots,s$, write $\tilde{f}_i = z_i - \tilde{h}_i$ with $\tilde{h}_i \in P$.
\begin{enumerate}
\item[(a)] The $K$-algebra homomorphism $\phi:\; P \longrightarrow K[Y]$
defined by $z_i \mapsto \tilde{h}_i$ for $i=1,\dots,s$ and $\phi(x_i)=x_i$
for $x_i \notin Z$ induces a $K$-algebra isomorphism
$$
\Phi:\; P/I \longrightarrow K[Y] / (I \cap K[Y])
$$
which is called a {\bf $Z$-separating re-embedding} of~$I$.

\item[(b)] Suppose that $I\subseteq \M= \langle x_1,\dots,x_s\rangle$
and $s=\dim_K(\Lin_{\M}(I))$. Then we have $\edim(P/I) = n-s$, i.e.,
the $Z$-separating re-embedding of~$I$ is an optimal re-embedding.
\end{enumerate}
\end{proposition}

In view of this proposition, the task of finding good, or even optimal, re-em\-bed\-dings
of a polynomial ideal can be divided into two steps:
\begin{enumerate}
\item[(1)] Find candidate tuples~$Z$ such that the ideal~$I$ has a reasonable
chance of containing a $Z$-separating tuple of polynomials.

\item[(2)] Given a candidate tuple~$Z$, check whether~$I$ does indeed contain
a $Z$-separating tuple.
\end{enumerate}

To perform the first step, for instance the following approaches have been proposed.

\begin{remark}\label{rem-findZ}
Let $I$ be a proper ideal in~$P$. Suppose we want to find a tuple~$Z$ of distinct
indeteminates in~$X$ which is a separating tuple for~$I$.
\begin{enumerate}
\item[(a)] If we can find a term ordering~$\sigma$ such that $Z \subseteq \LT_\sigma(I)$
then~$Z$ is a separating tuple of indeterminates for~$I$. For instance, we can take all
indeterminates contained in~$\LT_\sigma(I)$.

Unfortunately, checking all leading term ideals $\LT_\sigma(I)$ involves the computation
of the Gr\"obner fan of~$I$. This is known to be an extremely hard computation in general,
as the Gr\"obner fan tends to contain many elements. One improvement can be to compute
a {\it restricted Gr\"obner fan} (see~\cite{KLR2}).

\item[(b)] If the ideal $I$ is contained in~$\M = \langle x_1,\dots,x_n\rangle$, it turns out
that it is sufficient to compute the Gr\"obner fan of the linear part $\Lin_{\M}(I)$ of~$I$
(see~\cite{KLR3}, Prop.~2.6). This task can be reduced to finding the minors of a matrix
(see~\cite{KLR3}, Thm.~3.5 and Alg.~4.1). In general, this may still be quite demanding,
but for instance when $\Lin_{\M}(I)$ is a binomial linear ideal, it can be simplified further
(see~\cite{KLR3}, Alg.~5.7).
\end{enumerate}
\end{remark}

In view of this remark we concentrate in the following on step~(2) of the above procedure.
In other words, we assume that we are given a tuple~$Z$ of indeterminates and want to
check whether it is separating for~$I$.
Even if an ideal contains a $Z$-separating tuple of indeterminates, this may not
be apparent by looking at its generators (see for instance \cite{KLR3}, Example 2.4).
In general, checking whether a given tuple $Z$ of indeterminates is a separating
tuple of indeterminates for~$I$ may require us to compute a Gr\"obner basis
and thus be infeasible for large examples.

\bigbreak
%
%

\section{Checking Separating Tuples of Indeterminates}
\label{sec3-check}

As before, we let~$I$ be a proper ideal in $P=K[x_1,\dots,x_n]$,
we let $X=(x_1,\dots,x_n)$, and we let $Z=(z_1,\dots,z_s)$ be a tuple of distinct
indeterminates in~$X$.
Our goal in this section is to find fast methods
for determining whether $Z$ is a separating tuple of indeterminates for~$I$, i.e.,
whether~$I$ contains a $Z$-separating tuple of polynomials.
Since the calculation of a Gr\"obner basis of~$I$ with respect to an elimination
term ordering for~$Z$ is, in general, too expensive, we construct procedures which
are allowed to fail.
We start by introducing the following algorithm, which is a part of
Algorithm~\ref{alg-checkZ}, the main algorithm of this paper.

\begin{algo}{\bf (Linear Interreduction)}\label{alg-linInterred}\\
Let $Z=(z_1,\dots,z_s)$ be a tuple of distinct indeterminates in~$X$, and let $\sigma$
be the lexicographic term ordering on~$P$.
Moreover, let $g_1,\dots,g_r \in P$ such that
$\dim_K(\langle \Lin(g_1),\dots,\Lin(g_r) \rangle_K)=s$
and such that every term in the union of the supports
of $g_1,\dots,g_r$ is divisible by at least one of the indeterminates in~$Z$.
Consider the following sequence of instructions.

\begin{algorithmic}[1]
\State Compute $T = (z_1,\dots,z_s,t_1,\dots,t_m)$ where $\{t_1,\dots,t_m\}$ is the combined support
of $g_1,\dots,g_r$ without $z_1,\dots,z_s$, ordered such that $t_1 >_\sigma \dots >_\sigma t_m$. 

\State Form the coefficient matrix $M \in \Mat_{r,m+s}(K)$ of $g_1,\dots,g_r$ with respect to
the terms in $T$.

\State Compute the matrix $N \in \Mat_{r',r}(K)$ such that $N \cdot M \in \Mat_{r',m+s}(K)$ is in
reduced row echelon form and has no zero rows. \;

\State Compute $(g_1',\dots,g_{r'}')^{\operatorname{tr}} = N \cdot (g_1,\dots,g_r)^{\operatorname{tr}}$
and return $g_1',\dots,g_{r'}'$.
\end{algorithmic}

This is an algorithm which computes a $K$-basis $g_1',\dots,g_{r'}'$ of $\langle g_1,\dots,g_r \rangle_K$ 
such that the following conditions hold.
\begin{itemize}
\item[(a)] For $i=1,\dots,s$, the polynomial $g'_i$ is of the form $g'_i = z_i - h_i$, where
$h_i\in P_{\ge 2}$. 

\item[(b)] Let $u=r'-s$. For $j=1,\dots,u$, we get $g'_{s+j}=q_j \in P_{\ge 2}$ such that
the polynomials $q_1,\dots,q_u$ have distinct leading terms with respect to~$\sigma$.

\item[(c)] The polynomials $q_1,\dots,q_u$ form a $K$-basis of 
$\langle g_1,\dots,g_r \rangle_K \cap P_{\ge 2}$.

\item[(d)] For every $i\in \{1,\dots,s\}$, we have $z_i \in \langle g_1,\dots,g_r \rangle_K$ 
if and only if $g_i' = z_i$.

\end{itemize}
\end{algo}

\begin{proof}
Let $U = \langle g_1,\dots,g_r \rangle_K$.
First observe that $N \cdot M$ is the coefficient matrix of $g'_1,\dots,g'_{r'}$ with
respect to the terms in~$T$.
Since the rows of $N \cdot M$ form a $K$-basis of the row space of~$M$, the polynomials
$g'_1, \dots, g'_{r'}$ form a $K$-basis of~$U$ by the isomorphism between the row space
of $M$ and $U$.

To prove (a), we first note that $1 \notin \{t_1,\dots,t_m\}$ since $1$ is not divisible by
any indeterminate in~$Z$.
Let $V = \langle \Lin(g_1),\dots,\Lin(g_r) \rangle_K$. Since $\dim_K(V) = s$ and every
term in the combined support is divisible by an indeterminate from $Z$, we have
$V = \langle Z \rangle_K$. This shows that the upper left $s \times s$ submatrix of
$N \cdot M$ is the identity matrix and hence $\Lin(g_i') = z_i$ for $i=1,\dots,s$
and we can write $g_i' = z_i-h_i$ where $h_i \in P_{\geq 2}$ since $1 \notin \Supp(h_i)$.

To prove (b), we note that $\Lin(g'_{s+j}) = 0$ and $1 \notin \Supp(g'_{s+j})$ for $j=1,\dots,u$.
This yields $g_{s+j}' \in P_{\ge 2}$. The remaining claim follows from the observations
that the coefficient matrix $N \cdot M$ of
$g_1',\dots,g_{r'}'$ is in reduced row-echelon form and that the terms $t_1,\dots,t_m$
are sorted decreasingly with respect to $\sigma$.

For the proof of~(c), let $q \in U \cap P_{\ge 2}$. By~(a), there exist
$a_1,\dots,a_s, b_1,\dots,b_u \in K$ such that
$q = a_1(z_1-h_1) + \dots + a_s(z_s-h_s)+b_1q_1+\dots+b_uq_u$.
Then $0 = \Lin(q) = a_1z_1+\dots+a_sz_s$ implies $a_1=\dots=a_s = 0$. Hence
we get $q \in \langle q_1,\dots,q_u \rangle_K$.
This shows $\langle q_1,\dots,q_u \rangle_K = U \cap P_{\ge 2}$. By (b), the polynomials
$q_1,\dots,q_u$ are $K$-linearly independent, and therefore a $K$-basis of $U \cap P_{\ge 2}$.

Finally, we show~(d). Let $i \in \{1,\dots,s\}$ such that $z_i \in U$. Since $g_i' \in U$, we get
$h_i = z_i-g_i' \in U \cap P_{\ge 2}$. Then~(c) yields $h_i \in \langle q_1,\dots,q_u \rangle_K$.
Since $N \cdot M$ is in reduced row echelon form, the leading terms of $q_1,\dots,q_u$
do not appear in the support of~$h_i$. However, since these polynomials have distinct leading terms
according to~(b), every nonzero element $q \in \langle q_1,\dots,q_u \rangle_K$
satisfies $\LT_\sigma(q) = \LT_\sigma(q_j)$ for some~$j$. This shows $h_i = 0$,
and hence $g_i' = z_i$. For the other implication, note that if $g_i' = z_i$,
then clearly $z_i = g_i' \in U$.
\end{proof}

Now we are ready to introduce the main algorithm of this paper.
Its basic idea is to try to find a tuple of 
weights $W=(w_1,\dots,w_n)$ for the indeterminates in~$X$ such that
for every term ordering~$\sigma$ compatible with the grading given by~$W$ there exist
polynomials $f_i \in I$ with $\LT_\sigma(f_i)=z_i$ for $i=1,\dots,s$.

\begin{algo}{\bf (Checking Separating Indeterminates)}\label{alg-checkZ}\\
Let $g_1,\dots,g_r \in P$, let $I = \langle g_1,\dots,g_r\rangle$,
and let $Z=(z_1,\dots,z_s)$ be a tuple of distinct indeterminates in~$X$.
Consider the following sequence of instructions.

\begin{algorithmic}[1]
\State Set $w_i = 0$ for $i \in \{1,\dots,n\}$, $\delta = \max\{\deg(g_1),\dots,\deg(g_r)\}$,
and $d=1$.
\label{algstep-init}
\State Delete in~$g_1,\dots,g_r$ every monomial which is not divisible
by an indeterminate in~$Z$ and call the result $g_1,\dots,g_r$ again.
\label{algstep-delterms1}
\If{$\dim_K( \langle \Lin(g_1), \dots, \Lin(g_r) \rangle_K) < \# Z$} \label{algstep-if1-start}
   \State \Return ``\texttt{Fail}''. \label{algstep-if1-fail}
\EndIf \label{algstep-if1-end}

\Repeat \label{algstep-loop-start}
  \State Apply Algorithm \ref{alg-linInterred} to $g_1,\dots,g_r$ and $Z$
  and call the result $g_1,\dots,g_r$ again. \label{algstep-interred}
  \State Set $\widetilde{Z} = Z \cap \{g_1,\dots,g_r\}$. \label{algstep-tildeZ}
  \If{$\widetilde{Z} = \emptyset$} \label{algstep-ifloop-start}
    \State \Return ``\texttt{Fail}''. \label{algstep-ifloop-fail}
  \EndIf \label{algstep-ifloop-end}
  \ForAll{$z_i \in \widetilde{Z}$} \label{algstep-forzi-start}
    \State Let $k \in \{1,\dots,n\}$ with $x_k = z_i$ and set $w_k = d$. \label{algstep-setweight}
    \State Remove $z_i$ from $Z$. \label{algstep-removezi}
  \EndFor \label{algstep-forzi-end}
  \State Delete in $g_1,\dots,g_r$ every monomial which is not divisible by an
  indeterminate \hbox{\qquad} in~$Z$ and call the result $g_1,\dots,g_r$ again. \label{algstep-delterms2}
  \State Replace $d$ by $\delta\cdot d+1$. \label{algstep-setd}
\Until{$Z$ is empty.} \label{algstep-loop-end}
\State\Return $W=(w_1,\dots,w_n)$. \label{algstep-return}
\end{algorithmic}

\noindent This is an algorithm which returns either ``\texttt{Fail}'' or
a tuple of non-negative weights $W=(w_1,\dots,w_n) \in \mathbb{N}^n$.
If it returns a tuple $W$, then there are $f_1,\dots,f_s \in \langle g_1,\dots,g_r \rangle_K$
such that $(f_1,\dots,f_s)$ is a $Z$-separating tuple of polynomials for~$I$,
and such that every term ordering~$\sigma$ which is compatible with the grading
given by~$W$ satisfies $\LT_\sigma(f_i)=z_i$ for $i=1,\dots,s$.
\end{algo}

\begin{proof}
For the finiteness of the algorithm, it suffices to verify that
the loop in Steps~\ref{algstep-loop-start}-\ref{algstep-loop-end} terminates.
In every iteration of the loop, the algorithm either stops in Step~\ref{algstep-ifloop-fail}
or the number of elements in $Z$
decreases in Step~\ref{algstep-removezi}.
This yields that the loop terminates after finitely many iterations.
\smallskip

Now suppose that we are at the start of an iteration of the loop in Steps~\ref{algstep-loop-start}-\ref{algstep-loop-end}.
Let $U$ be the $K$-vector space generated by the input polynomials and let $T$
be the set of all terms deleted in Step \ref{algstep-delterms1} or in an execution of Step \ref{algstep-delterms2} in a
previous iteration of the loop. We show that the following properties hold
at the start of every iteration of the loop.
\begin{enumerate}
\item[(1)]
The polynomials $g_1,\dots,g_r$ are a valid input for Algorithm \ref{alg-linInterred}, i.e.,
   \begin{itemize}
   \item[(1.1)]
   every term in the support of $g_1,\dots,g_r$ is divisible by an indeterminate from $Z$,
   
   \item[(1.2)] $\dim_K(\langle \Lin(g_1),\dots,\Lin(g_r) \rangle_K) = \# Z$.
   \end{itemize}

\item[(2)]
We have $g_1,\dots,g_r \in U + \langle T \rangle_K$.

\item[(3)]
All terms in $T$ have $(w_1,\dots,w_n)$-weight smaller than $d$.
\end{enumerate}

To prove (1.1), we note that every monomial not divisible by an indeterminate from $Z$ was either
deleted in Step \ref{algstep-delterms1} or in Step \ref{algstep-delterms2} in an earlier iteration of the loop, so the claim holds. 

Next we show (1.2). In the first iteration, passing the condition in Steps~\ref{algstep-if1-start}-\ref{algstep-if1-end} guarantees the claim.
Assume that it holds at the start of an iteration of the loop that is not the
last iteration. Then, after Step~\ref{algstep-interred},
Algorithm~\ref{alg-linInterred}.a yields $h_1,\dots,h_s,q_1,\dots,q_r \in P_{\ge 2}$
such that $(g_1,\dots,g_r) = (z_1-h_1,\dots,z_s-h_s,q_1,\dots,q_r)$.
Without loss of generality assume that $\widetilde{Z} = \{z_1,\dots,z_k\}$ for $k < s$.
Then, after the execution of Step \ref{algstep-delterms2} and therefore at the start of the next loop, we have
$(g_1,\dots,g_r) = (0,\dots,0,z_{k+1}-\tilde{h}_{k+1},\dots,z_s-\tilde{h}_s,\tilde{q}_1,\dots,\tilde{q}_r)$
where $\tilde{h}_{k+1},\dots,\tilde{h}_s,\tilde{q}_1,\dots,\tilde{q}_r \in P_{\geq 2}$. 
Consequently, $\langle \Lin(g_1),\dots,\Lin(g_r) \rangle_K = \langle z_{k+1},\dots,z_s \rangle_K$
has dimension $s-k = \#Z$.

To show~(2), we observe that it holds after the execution of Step~\ref{algstep-delterms1}, and thus at the start of the first
iteration of the loop. Moreover, the property is preserved by $K$-linear interreductions and
monomial deletions, which are the only operations performed on $g_1,\dots,g_r$ during the loop
in Steps \ref{algstep-loop-start}-\ref{algstep-loop-end}. This implies~(2).

To prove~(3), notice first that, if a term $t$ is added to $T$ at any point of the algorithm, then
every later execution of Step~\ref{algstep-setweight} only changes the weight of an indeterminate in~$Z$
which does not divide~$t$. Hence the weight of~$t$ does not change anymore.
Moreover, if $t \in T$ is deleted in Step \ref{algstep-delterms1}, then its weight is~$0$ which is smaller than $d = 1$.
If it is deleted in Step \ref{algstep-delterms2}, then every indeterminate dividing $t$ has weight $0$ or has
been assigned a weight smaller than or equal to $d$. Hence $\deg(t) \leq \delta$ implies that
the weight of $t$ is smaller or equal to $\delta \cdot d < \delta \cdot d + 1$. This shows that,
at the start of the next iteration of the loop, the weight of~$t$ is smaller than~$d$.
\smallskip

Now assume that the algorithm successfully returned a weight tuple~$W$.
Then claim~(2) yields that, for every $z_i \in Z$, there was an iteration of the loop in Steps \ref{algstep-loop-start}-\ref{algstep-loop-end}
where we had $z_i \in U + \langle T \rangle_K$ in Step~\ref{algstep-tildeZ}. Thus we can write
$z_i = f_i+p_i$, where $f_i \in U$ and $p_i \in \langle T \rangle_K$.
By claim~(3), every term in $T$ has $W$-weight smaller than~$d$, while~$z_i$ has $W$-weight
$d$. Hence every term ordering~$\sigma$ compatible with the grading given by~$W$ satisfies
$\LT_\sigma(f_i) = \LT_\sigma(z_i-p_i) = z_i$. Notice that $f_i \in U = \langle g_1,\dots,g_r \rangle_K$.
\end{proof}

The next proposition characterizes when Algorithm \ref{alg-checkZ} is successful.

\begin{proposition}\label{prop-checkZ-success}
Let $g_1,\dots,g_r \in P$, let $I = \langle g_1,\dots,g_r \rangle$,
and let $Z = (z_1,\dots,z_s)$ be a tuple of distinct indeterminates in $X$.
Then the following conditions are equivalent.
\begin{itemize}
\item[(a)] Algorithm~\ref{alg-checkZ} successfully computes a weight tuple~$W$.

\item[(b)] There exists a $Z$-separating tuple $(f_1,\dots,f_s)$ of polynomials in~$I$
such that $f_1,\dots,f_s \in \langle g_1,\dots,g_r \rangle_K$.

\end{itemize}
\end{proposition}

\begin{proof}
The correctness of Algorithm \ref{alg-checkZ} shows that (a) implies~(b).
Thus we only need to show that (b) implies~(a).
Assume that there are $f_1,\dots,f_s \in \langle g_1,\dots,g_r \rangle_K$
such that $(f_1,\dots,f_s)$ is a $Z$-separating tuple of polynomials in~$I$.
Let $\tau$ be a term ordering on $P$ such that $\LT_\tau(f_i)=z_i$ for $i=1,\dots,s$.

First we show that the algorithm does not terminate in Step~\ref{algstep-if1-fail}.
For this, let $f_1',\dots,f_s'$ be the result of deleting in $f_1,\dots,f_s$ every
term which is not divisible by an indeterminate from $Z$.
Similarly, for $i=1,\dots,r$, let $g'_i$ be obtained from~$g_i$ in the same way.
Then $f_1',\dots,f_s' \in \langle g'_1,\dots,g'_r \rangle_K$.
This shows $\Lin(f_1'),\dots,\Lin(f_s') \in \langle \Lin(g_1),\dots,\Lin(g_r) \rangle_K$.
Here $\Lin(f_1'),\dots,\Lin(f_s')$ are $K$-linearly independent, since
$\LT_\tau(f_i') = z_i$ yields $\LT_\tau(\Lin(f_i')) = z_i$ for $i=1,\dots,s$.
Hence the condition in Step~\ref{algstep-if1-start} is false. 
\smallskip

It remains to show that the algorithm does not fail in Step~\ref{algstep-ifloop-fail}.
For this, let $\hat{g}_1, \dots, \hat{g}_{\hat{r}}$ be the original input polynomials of the algorithm,
let $\widehat{U} = \langle \hat{g}_1, \dots, \hat{g}_{\hat{r}} \rangle_K$, 
and write $\widehat{S} = \bigcup_{i=1}^{\hat{r}} \Supp(\hat{g}_i)$.
At any state of the algorithm, let $U = \langle g_1,\dots,g_r \rangle$,
let $S = \bigcup_{i=1}^r \Supp(g_i)$,
let $T = \widehat{S} \setminus S$ be the set of deleted terms,
and let $\phi : \langle S' \rangle_K \to \langle S \rangle_K$
be the map representing the deletions, i.e., the $K$-linear map
defined by $\phi(t) = t$ for $t \in S$ and $\phi(t) = 0$ for $t \in T$.

Let us show that $U = \phi(\widehat{U})$ at every step of the algorithm. 
First notice that after the execution of Step~\ref{algstep-if1-start}, we have $r=\hat{r}$ and
$g_i = \phi(\hat{g}_i)$ for $i=1,\dots,r$, so
$U = \langle g_1,\dots,g_r \rangle_K = \langle \phi(\hat{g}_1), \dots, \phi(\hat{g}_{\hat{r}}) \rangle_K 
= \varphi(\widehat{U})$.
This still holds at the start of the first iteration of the loop in Steps \ref{algstep-loop-start}-\ref{algstep-loop-end}.
Since $U$ is invariant under Step~\ref{algstep-interred} in view of Algorithm \ref{alg-linInterred}, it is only left
to show that the claim still holds after the execution of Step~\ref{algstep-delterms2}.
For this, assume that we are after an execution of Step~\ref{algstep-delterms2}. Write $\bar S$,
$\bar \phi$, and $\bar U$ for the values of~$S$, $\phi$, and~$U$ before
the execution of Step~\ref{algstep-delterms2}, respectively.
Write $\psi : \langle \bar S \rangle_K \to \langle S \rangle_K$ for the map
representing the term deletion in this single execution of Step~\ref{algstep-delterms2}, i.e.,
the $K$-linear map with $\psi(t) = t$ for $t \in S$ and $\psi(t) = 0$ for $t \in \bar S \setminus S$.
Then $\phi = \psi \circ \bar{\phi}$ and $U = \psi(\bar U) = \psi(\bar{\phi}(\widehat{U})) = \phi(\widehat{U})$.
Therefore $U = \phi(\widehat{U})$ holds after the execution of Step~\ref{algstep-delterms2}, and hence at the start
of the next iteration of the loop. 
\smallskip

Now we show that $\widetilde{Z} \ne \emptyset$ in every execution of Step~\ref{algstep-tildeZ}.
Assume that we are at the start of an iteration of the loop in 
Steps \ref{algstep-loop-start}-\ref{algstep-loop-end} and let
$z_i \in Z$ be the minimal element of~$Z$ with respect to~$\tau$.
Write $f_i = az_i+p_i$ with $a \in K \setminus \{0\}$ and $p_i \in P$ such that $az_i$
is the leading monomial of $f_i$. Then none
of the indeterminates in~$Z$ divides any of the terms in~$p_i$, since any term
divisible by an element of~$Z$ is larger than or equal to~$z_i$ with respect to~$\tau$.
This shows $\Supp(p_i) \cap S = \emptyset$, which in turn implies $\phi(p_i) = 0$, and hence
$z_i = \phi(f_i) \in \phi(\widehat{U}) = U = \langle g_1,\dots,g_r \rangle_K$.
Finally, looking at Algorithm~\ref{alg-linInterred}.d, we note that after the execution 
of Step~\ref{algstep-interred} 
we have $g_j = z_i$ for some index $j \in \{1,\dots,r\}$, and hence $g_j \in \widetilde{Z}$ 
in Step~\ref{algstep-tildeZ}.
\end{proof}

The following example illustrates Algorithm~\ref{alg-checkZ} at work.
It will also be reconsidered in the next sections.

\begin{example}\label{ex-checkZsepI}
Let $P = \QQ[x_1,\dots,x_{11}]$, and let~$I$ be the ideal in~$P$ generated by the
polynomials $g_1,\dots,g_9 \in P$ given by
\begin{align*}
  g_1 &= x_{1} x_{4} x_{8} x_{11} +x_{5} x_{6} x_{8} +x_{5} x_{6} x_{10}
  	+x_{3} x_{6} +x_{7} x_{8} +x_{1} +x_{4}, \\
  g_2 &=  -x_{1}^2 x_{3}^4  -x_{1} x_{2} x_{3} x_{6} x_{8} x_{11}
  	+x_{1}^3 x_{10} +x_{3} +x_{7} +1, \\
  g_3 &= x_{1} x_{2} x_{3}^2 x_{8}^2  +x_{6} x_{7}^2 x_{8} +x_{1} x_{3} x_{8}^2
  	+x_{7} x_{8} +x_{7} x_{10} +x_{2} +x_{5}, \\
  g_4 &=  -x_{1} x_{7} x_{9} +x_{3} x_{11} +x_{9}, \\
  g_5 &= x_{1} x_{4} x_{8} x_{11} +x_{1}^2 x_{6} +x_{5} x_{7}, \\
  g_6 &= x_{8} x_{10}^2  +x_{7} x_{9}, \\
  g_7 &= x_{1}^2 x_{3}^4  +x_{2}^2 x_{6}^4  +x_{1} x_{2} x_{3}^2 x_{8}^2
  	+x_{1} x_{2} x_{6}^2 x_{10}^2  +x_{1} x_{2} x_{3} x_{6} x_{8} x_{11}, \\
  g_8 &= x_{2}^6  +x_{1} x_{6} x_{8} -x_{3} x_{6}, \\
  g_9 &= x_{1}^6  +x_{1}.
\end{align*}
We apply Algorithm \ref{alg-checkZ} to check whether $Z = (x_4,x_5,x_7)$
is a separating tuple of indeterminates for $I$.
\begin{itemize}
\item[(1)]
In Step~\ref{algstep-init}, we set $w_1=\dots=w_{11}=0$, $\delta =6$, and $d = 1$.

\item[(2)] In Step~\ref{algstep-delterms1}, we delete all monomials in $g_1,\dots,g_9$ not divisible by
$x_4$, $x_5$, or $x_7$, and obtain
\begin{align*}
  g_1 &= x_{1}x_{4}x_{8}x_{11}  + x_{5}x_{6}x_{8} + x_{5}x_{6}x_{10}
    	+ x_{7}x_{8} + x_{4}, \quad
  g_2 = x_{7}, \\
  g_3 &= x_{6}x_{7}x_{8} + x_{7}x_{8} + x_{7}x_{10} + x_{5}, \quad
  g_4 = -x_{1}x_{7}x_{9}, \\
  g_5 &= x_{1}x_{4}x_{8}x_{11} + x_{5}x_{7}, \quad
  g_6 = x_{7}x_{9}, \quad g_7=g_8=g_9=0.
\end{align*}
We have $\dim_\QQ(\langle \Lin(g_i) \mid i=1,\dots,9 \rangle_\QQ)=3=\#Z$,
so the algorithm does not stop in Step~\ref{algstep-delterms1}.

\item[(3)]
The result of applying Algorithm~\ref{alg-linInterred} in Step~\ref{algstep-interred} is
\begin{align*}
  g_1 &= x_{4} - x_{5}x_{7}  + x_{5}x_{6}x_{8} + x_{5}x_{6}x_{10} + x_{7}x_{8}, \\
  g_2 &= x_{5} + x_{6}x_{7}x_{8} + x_{7}x_{8} + x_{7}x_{10}, \\
  g_3 &= x_{7}, \quad
  g_4 = x_{1}x_{4}x_{8}x_{11} + x_{5}x_{7}, \\
  g_5 &= -x_{1}x_{7}x_{9},\quad
  g_6 = x_{7}x_{9}.
\end{align*}

\item[(4)]
In Step \ref{algstep-tildeZ} we obtain $\widetilde{Z} = \{x_7\}$.

\item[(5)]
Since $\widetilde{Z}\ne \emptyset$, Step \ref{algstep-setweight} sets $w_7 = 1$ and Step \ref{algstep-removezi} sets $Z = (x_4,x_5)$.
Deleting all monomials not divisible by $x_4$ or $x_5$ in $g_1,\dots,g_6$ in Step~\ref{algstep-delterms2} yields
\begin{align*}
  g_1 &= x_{4} - x_{5}x_{7}  + x_{5}x_{6}x_{8} + x_{5}x_{6}x_{10},\quad
  g_2 = x_{5},\quad g_3=0, \\
  g_4 &= x_{1}x_{4}x_{8}x_{11} + x_{5}x_{7},\quad
  g_5 = g_6 = 0.
\end{align*}
Step \ref{algstep-setd} sets $d$ to~$7$, and the loop repeats.

\item[(6)]
Algorithm \ref{alg-linInterred} returns
\begin{align*}
  g_1 &= x_{4} - x_{5}x_{7}  + x_{5}x_{6}x_{8} + x_{5}x_{6}x_{10}, \\
  g_2 &= x_{5}, \quad g_3 = x_{1}x_{4}x_{8}x_{11} + x_{5}x_{7}.
\end{align*}
Step \ref{algstep-tildeZ} finds $\widetilde{Z} =\{x_5\} \ne \emptyset$, so Step \ref{algstep-setweight} sets $w_5 = 7$, and
Step \ref{algstep-removezi} sets $Z = (x_4)$.
The deletions in Step \ref{algstep-delterms2} yield $g_1 = x_{4}$, $g_2 = 0$, and $g_3 = x_{1}x_{4}x_{8}x_{11}$.
Step \ref{algstep-setd} sets $d = 43$ and the loop repeats.

\item[(7)]
Algorithm \ref{alg-linInterred} returns $g_1 = x_4$ and $g_2 = x_{1}x_{4}x_{8}x_{11}$.
We have $\widetilde{Z} =\{x_4\}$ in Step \ref{algstep-tildeZ},
so Step \ref{algstep-setweight} sets $w_4 = 43$ and Step \ref{algstep-removezi} sets $Z = \emptyset$. After executing Steps~\ref{algstep-delterms2} and~\ref{algstep-setd},
the loop terminates and the algorithm returns $W = (0, 0, 0, 43, 7, 0, 1, 0, 0, 0, 0)$.
\end{itemize}
Hence there exists an $(x_4,x_5,x_7)$-separating tuple of polynomials in~$I$.

It is worth noting that the computation of a $\tau$-Gr\"obner basis of~$I$ 
from $g_1,\dots,g_9$ for an elimination ordering~$\tau$ for~$Z$ is quite hard and takes
many hours using computer algebra systems such as~\texttt{CoCoA} (see \cite{CoCoA}),
\texttt{msolve} (see \cite{msolve}), or \texttt{SageMath} (see \cite{sage}).
Algorithm~\ref{alg-checkZ}, however, is done in less than a second.
\end{example}

The next example shows an application of Algorithm~\ref{alg-checkZ} to
an ideal whose linear part is a binomial linear ideal
which was studied in \cite[Section~5]{KLR3}. This again shows that
our method can work efficiently with rather large examples where any Gr\"obner basis
calculation would be infeasible.

\begin{example}
In the polynomial ring $P=\QQ[x_{1},\dots,x_{84}]$, consider the ideal $I$ generated
by nonzero entries of the pairwise commutators of the following three matrices
\begin{center}
$\left( \begin{array}{ccccccc}
   0 & x_{3} & x_{5} &  0 &  0 &  0 & x_{12} \\
   0 & x_{15} & x_{17} &  0 &  0 &  0 & x_{24} \\
   0 & x_{27} & x_{29} &  0 &  0 &  0 & x_{36} \\
  1 & x_{39} & x_{41} &  0 &  0 &  0 & x_{48} \\
   0 & x_{51} & x_{53} & 1 &  0 &  0 & x_{60} \\
   0 & x_{63} & x_{65} &  0 & 1 &  0 & x_{72} \\
   0 & x_{75} & x_{77} &  0 &  0 & 1 & x_{84}
\end{array}\right),\qquad
\left( \begin{array}{ccccccc}
   0 & x_{2} & x_{4} & x_{5} & x_{7} & x_{9} & x_{11} \\
   0 & x_{14} & x_{16} & x_{17} & x_{19} & x_{21} & x_{23} \\
  1 & x_{26} & x_{28} & x_{29} & x_{31} & x_{33} & x_{35} \\
   0 & x_{38} & x_{40} & x_{41} & x_{43} & x_{45} & x_{47} \\
   0 & x_{50} & x_{52} & x_{53} & x_{55} & x_{57} & x_{59} \\
   0 & x_{62} & x_{64} & x_{65} & x_{67} & x_{69} & x_{71} \\
   0 & x_{74} & x_{76} & x_{77} & x_{79} & x_{81} & x_{83}
\end{array}\right),$
\\
$\left( \begin{array}{ccccccc}
   0 & x_{1} & x_{2} & x_{3} & x_{6} & x_{8} & x_{10} \\
  1 & x_{13} & x_{14} & x_{15} & x_{18} & x_{20} & x_{22} \\
   0 & x_{25} & x_{26} & x_{27} & x_{30} & x_{32} & x_{34} \\
   0 & x_{37} & x_{38} & x_{39} & x_{42} & x_{44} & x_{46} \\
   0 & x_{49} & x_{50} & x_{51} & x_{54} & x_{56} & x_{58} \\
   0 & x_{61} & x_{62} & x_{63} & x_{66} & x_{68} & x_{70} \\
   0 & x_{73} & x_{74} & x_{75} & x_{78} & x_{80} & x_{82}
\end{array}\right).$
\end{center}
Then a system of generators of $I$ consists of 126 quadratic polynomials
without constant, and so $I\subseteq \M$.
The linear part of~$I$ has dimension $57$ and generates a binomial linear ideal.
Using the method mentioned in Remark~\ref{rem-findZ}.b,
we find that the following tuple $Z$ of 57 indeterminates
\begin{align*}
(& x_{1}, x_{2}, x_{3}, x_{4}, x_{5}, x_{6}, x_{7}, x_{8}, x_{9}, x_{10},
x_{11}, x_{12}, x_{18}, x_{19}, x_{20}, x_{21}, x_{22}, x_{23}, x_{30}, \\
& x_{31}, x_{32}, x_{33}, x_{34}, x_{35}, x_{37}, x_{38}, x_{40}, x_{42},
x_{43}, x_{44}, x_{45}, x_{46}, x_{47}, x_{49}, x_{50}, x_{52}, x_{54}, x_{55},\\
& x_{56}, x_{57}, x_{58}, x_{59}, x_{61}, x_{62}, x_{64}, x_{66}, x_{67},
x_{68}, x_{69}, x_{70}, x_{71}, x_{78}, x_{79}, x_{80}, x_{81}, x_{82}, x_{83})
\end{align*}
is a possible candidate for a separating tuple of indeterminates for~$I$.
When we apply Algorithm~\ref{alg-checkZ}, it shows
that~$Z$ is indeed a separating tuple of indeterminates for~$I$ and returns
the weight tuple
\begin{align*}
(& 127,  127,  15,  127,  15,  31,  31,  31,  31,  31,  31,  15,
1,  1,  3, 3,  7,  7,  1, \\
& 1,  3,  3,  7,  7, 63,  63,  63,  15,  15,  31,  31,  31,  31, 31,  31,  31,  1,\\
& 1,  15,  15,  31,  31, 15,  15,  15,  1,  1,  3,  3,  15,  15,
1,  1,  3,  3,  7,  7).
\end{align*}
Notice that we also have $\# Z=57 = \dim_K(\Lin_\M(I))$.
Thus the $Z$-separating re-embedding of~$I$ is an optimal re-embedding
by~\cite[Corollary 2.8]{KLR3}.

Moreover, the closed subscheme of $\AA^{84}_\QQ$ defined by~$I$
is the border basis scheme associated to a certain order ideal (see \cite[Example 8]{BC}).
Using techniques relying on Pommaret bases, it is shown in~\cite{BC}
that this scheme can be embedded into $\AA^{42}_\QQ$, while our method
yields an optimal embedding into~$\AA^{27}_\QQ$.
\end{example}

\bigbreak
%
%

\section{Optimized Checking of Separating Indeterminates}
\label{sec4-checkopt}

In this section we discuss an optimization of Algorithm~\ref{alg-checkZ}.
Instead of simply interreducing the polynomials $g_1,\dots,g_r$ linearly using
Algorithm~\ref{alg-linInterred} in Step~\ref{algstep-interred}, we first extend them with some
polynomials contained in $P_{\ge 2}$.
This allows Algorithm~\ref{alg-OptimizedCheckZ} to succeed in cases where
Algorithm~\ref{alg-checkZ} fails.

\begin{algo}\label{alg-OptimizedCheckZ}
Let $I = \langle g_1,\dots,g_r\rangle$ be a proper ideal in~$P$, where
$g_1,\dots,g_r \in P \setminus \{0\}$,
and let $Z=(z_1,\dots,z_s)$ be a tuple of distinct indeterminates in~$X$.
Consider the following modifications of Algorithm~\ref{alg-checkZ}. 
\smallskip

\noindent (I) Replace Step~\ref{algstep-interred} by the following sequence of steps.

\begin{algorithmic}
\State \hspace{-0.4cm}{\footnotesize \ref{algstep-interred}a:} Form the set $H=\{x_i g_j \mid x_i \in X \setminus Z \text{ and } j \in \{1,\dots,r\}\}$.

\State \hspace{-0.4cm}{\footnotesize \ref{algstep-interred}b:} Let $\tau$ be a degree-compatible term ordering.
Compute a $K$-basis $\{h_1,\dots,h_m\}$ of the space
$\langle H \rangle_K$ such that $\LT_\tau(h_i) \ne \LT_\tau(h_j)$ for $i \ne j$.

\State \hspace{-0.4cm}{\footnotesize \ref{algstep-interred}c:} Compute the intersection $\{q_1,\dots,q_u\} = \{h_1,\dots,h_m\} \cap P_{\le\delta}$.

\State \hspace{-0.4cm}{\footnotesize \ref{algstep-interred}d:} Apply Algorithm \ref{alg-linInterred} to $Z$ and $g_1,\dots,g_r,q_1,\dots,q_u$.
Call the result $g_1,\dots,g_r$ again.
\end{algorithmic}
\smallskip

\noindent (II) Moreover, replace Step \ref{algstep-setd} by the following step:

\begin{algorithmic}
\State \hspace{-0.4cm}{\footnotesize \ref{algstep-setd}:} Replace $d$ by $2 \delta d + 1$.
\end{algorithmic}

\noindent 
Then the resulting sequence of instructions is an algorithm
which returns either ``\texttt{Fail}'' or
a tuple of non-negative weights $W=(w_1,\dots,w_n) \in \mathbb{N}^n$.
If it returns a tuple $W$, then there is a $Z$-separating tuple of polynomials
$(f_1,\dots,f_s)$ for~$I$
such that every term ordering~$\sigma$ which is compatible with the grading
given by~$W$ satisfies $\LT_\sigma(f_i)=z_i$ for $i=1,\dots,s$.
\end{algo}

\begin{proof}
The finiteness of the algorithm follows in the same way as in the proof of
Algorithm~\ref{alg-checkZ}, so we only show its correctness whenever it returns
a weight tuple $W=(w_1,\dots,w_s)$.
Also observe that Step \ref{algstep-interred}b can be executed e.g.~by executing Algorithm~\ref{alg-linInterred}
on the input $H$ and the empty tuple.
\smallskip

Let $\delta = \max\{\deg(g_1),\dots,\deg(g_r)\}$ as in Step \ref{algstep-init} of the algorithm.
For $k\in\NN$, let $U_k= I_{\le \delta+k}$ be the $K$-vector space of all
polynomials in~$I$ of degree $\le \delta+k$.
At any state of the algorithm, let $W = (w_1,\dots,w_n)$, and for every $\ell\ge 0$ 
let $V_\ell$ be the $K$-vector space generated by all terms in
$K[X \setminus Z]_{\le \delta+\ell}$ having $W$-weight $\le\ell$.
Notice that in every state of the algorithm, except between Steps \ref{algstep-setweight} and \ref{algstep-removezi},
every indeterminate in $X \setminus Z$
has already been assigned a weight. So, if a term is in $V_\ell$, then its $W$-weight
does not change in a later step of the algorithm, even if some weights $w_k$ are adjusted in
a later execution of Step~\ref{algstep-setweight}.
Moreover, both sequences $(U_k)_{k \in \NN}$ and $(V_\ell)_{\ell \in \NN}$ are
increasing sequences, i.e., $U_k \subseteq U_{k+1}$ and $V_\ell \subset V_{\ell+1}$
for all $k,\ell \in \NN$.

Next we prove the following five claims.

\begin{itemize}
\item[(1)]
At the start of the first iteration of the loop in Steps \ref{algstep-loop-start}-\ref{algstep-loop-end}, we have
$g_1,\dots,g_r \in U_0+V_0$.

\item[(2)]
If $g_1,\dots,g_r \in U_m+V_{m'}$
for some $m,m' \in \NN$ at the start of an iteration of the loop,
then we have $g_1,\dots,g_r \in U_{m+1}+V_{m'+(d-1)/2}$
after the execution of Steps \ref{algstep-interred}a-\ref{algstep-interred}d.

\item[(3)]
If $g_1,\dots,g_r \in U_{m}+V_{m'}$
for some $m,m' \in \NN$ after an execution of Steps \ref{algstep-interred}a-\ref{algstep-interred}d during an iteration of the loop,
then we have $g_1,\dots,g_r \in U_{m}+V_{m'}+V_{\delta d}$
after the execution of Step~\ref{algstep-delterms2}.

\item[(4)]
For every $k \in \NN$ such that there exists a $k$-th iteration of the loop in Steps \ref{algstep-loop-start}-\ref{algstep-loop-end},
we have $g_1,\dots,g_r \in U_{k-1}+V_{(d-1)/2}$ at the start of the $k$-th iteration.

\item[(5)]
For every $k \in \NN$ such that there exists a $k$-th iteration of the loop in Steps \ref{algstep-loop-start}-\ref{algstep-loop-end},
we have $g_1,\dots,g_r \in U_{k}+V_{d-1}$ after the execution of Steps \ref{algstep-interred}a-\ref{algstep-interred}d in the
$k$-th iteration of the loop.
\end{itemize}


To show (1), observe that all terms deleted in Step~\ref{algstep-delterms1} are elements of
$K[X \setminus Z]_{\leq \delta}$ and have $W$-weight~$0$.
This shows $g_1,\dots,g_r \in U_0+V_0$ at the start of the first iteration of the loop
in Steps \ref{algstep-loop-start}-\ref{algstep-loop-end}.
\smallskip


For a proof of~(2), assume that we are after the execution of Step~\ref{algstep-interred}c, and let
$U = \langle H \rangle_K \cap P_{\leq \delta}$. We show the following.

\noindent (2.1) The polynomials $q_1,\dots,q_u$ computed in Step~\ref{algstep-interred}c form a $K$-basis of~$U$.

\noindent (2.2) If $g_1,\dots,g_r \in U_m+V_{m'}$ for some $m,m' \in \NN$ before the execution of
Step~\ref{algstep-interred}a, then we have $U \subseteq U_{m+1} + V_{m'+(d-1)/2}$ after the execution of Step~\ref{algstep-interred}c.
\smallskip

To verify (2.1), we note that $q_1,\dots,q_u \in U$.
Let $f \in U$ and suppose without loss of generality
that $\LT_\tau(h_1) <_\tau \dots <_\tau \LT_\tau(h_m)$ in Step~\ref{algstep-interred}b.
Since $\tau$ is a degree-compatible term ordering, we have
$\{q_1,\dots,q_u\} = \{h_1,\dots,h_u\}$.
As $f \in \langle h_1,\dots,h_m \rangle_K$ and the leading terms of $h_1,\dots,h_m$
are pairwise distinct, there is an index $j \in \{1,\dots,m\}$ such that $\LT_\tau(f) = \LT_\tau(h_j)$.
It follows that $f\in \langle h_1,\dots,h_j \rangle_K$,
because, for all $i > j$, we have $\deg(h_i) > \deg(f)$, and hence $\LT_\tau(h_i) >_\tau \LT_\tau(f)$.
Moreover, $f \in P_{\le\delta}$ implies $j \le u$. This yields
$f\in \langle h_1,\dots,h_u \rangle_K= \langle q_1,\dots,q_u \rangle_K$,
and (2.1) follows.

For a proof of (2.2), it suffices to show that in Step~\ref{algstep-interred}a 
we have $x_ig_j \in  U_{m+1} +V_{m'+(d-1)/2}$ for 
$x_i \in X \setminus Z$ and $j \in \{1,\dots,r\}$.
By the assumption, $g_j \in U_{m}+V_{m'}$, so there are polynomials $u_j \in U_m$ and $v_j \in V_{m'}$
with $g_j = u_j+v_j$.
Every term in the support of $v_j$ has $W$-weight $\le m'$.
Moreover, the $W$-weight of $x_i$ was assigned in a previous iteration of the loop, so
it is at most $(d-1)/2\delta$, and hence smaller than $(d-1)/2$.
Altogether, the $W$-weight of every term in the support of $x_iv_j$
is at most $m'+(d-1)/2$ and we conclude that $x_iv_j \in V_{m' + (d-1)/2}$.
It follows that
$$
x_ig_j \;=\; x_iu_j + x_iv_j \;\in\; x_i \cdot U_m + V_{m'+(d-1)/2} \;\subseteq\; U_{m+1} + V_{m'+(d-1)/2}.
$$
Hence we get $U \subseteq U_{m+1}+V_{m'+(d-1)/2}$.
Now observe that after Step~\ref{algstep-interred}c we have
$$
g_1,\dots,g_r,q_1,\dots,q_u \;\in\; (U_m+V_{m'}) \cup (U_{m+1} + V_{m'+(d-1)/2})
\;=\; U_{m+1}+V_{m'+(d-1)/2},
$$
so after the execution of Step \ref{algstep-interred}d, we have $g_1,\dots,g_r \in U_{m+1}+V_{m'+(d-1)/2}$.
This shows Claim~(2).
\smallskip


To show (3), suppose that $g_1,\dots,g_r \in U_{m}+V_{m'}$ after an execution of Step~\ref{algstep-interred}d.
During the execution of Steps \ref{algstep-tildeZ}-\ref{algstep-forzi-end}, only $V_{m'}$ is extended. So, before the execution of
Step~\ref{algstep-delterms2}, we still have $g_1,\dots,g_r \in U_{m}+V_{m'}$.
Let $j \in \{1,\dots,r\}$, let
$c_1t_1, \dots, c_e t_e$ be the monomials deleted from $g_j$ in Step~\ref{algstep-delterms2},
where $c_1,\dots,c_e \in K$ and $t_1,\dots,t_e \in \TT^n$ are terms,
and let $g_j' = g_j - c_1 t_1 - \dots - c_e t_e$ be the resulting polynomial obtained from~$g_j$
after their deletion.
Notice that, for $i=1,\dots,e$, we have $t_i \in K[X \setminus Z]$
and $\deg(t_i) \le \delta$ by (2.1), so all $t_1,\dots,t_e$ have weight $\le \delta d$.
This shows $t_1,\dots,t_e \in V_{\delta d}$. Hence $g_j \in U_m+V_{m'}$ implies
$g_j' \in U_m+V_{m'}+V_{\delta d}$.
Consequently, after the execution of Step~\ref{algstep-delterms2}, we have $g_j \in U_m+V_{m'}+V_{\delta d}$
for all $i \in \{1,\dots,r\}$.
\smallskip


Next Claim (4) is shown inductively as follows. For $k=1$, it holds by~(1).
Moreover, let $k \in \NN$ such that there exist a $k$-th and a $(k+1)$-st iteration of the
loop, and suppose that $g_1,\dots,g_r \in U_{k-1}+V_{(d-1)/2}$ at the start of the
$k$-th iteration.
Then~(2) shows $g_1,\dots,g_r \in U_{k}+V_{(d-1)/2+(d-1)/2} = U_k+V_{d-1}$ after the
execution of Steps \ref{algstep-interred}a-\ref{algstep-interred}d of the iteration, and (3) implies
$$
g_1,\dots,g_r ~\in~ U_{k}+V_{d-1}+V_{\delta d} ~=~ U_{k}+V_{\delta d}
$$
after the execution of Step~\ref{algstep-delterms2}.
Therefore, after the adjustment of~$d$ in Step~\ref{algstep-setd} and at the start of the $(k+1)$-st loop,
we have $g_1,\dots,g_r \in U_k+V_{(d-1)/2}$. This shows~(4).
\smallskip


To prove (5), we note that at the start of the $k$-th iteration of the loop in Steps~\ref{algstep-loop-start}-\ref{algstep-loop-end},
we have $g_1,\dots,g_r \in U_{k-1}+V_{(d-1)/2}$ by~(4). Now (2) says that, after the execution of
Steps \ref{algstep-interred}a-\ref{algstep-interred}d, we have $g_1,\dots,g_r \in U_{k}+V_{(d-1)/2+(d-1)/2} = U_k+V_{d-1}$.
\smallskip


Finally, we combine everything and finish the correctness proof of the algorithm.
For every $z_i \in Z$, we eventually find $g_j = z_i$ for some index $j\in\{1,\dots,r\}$,
where we are in Steps~\ref{algstep-forzi-start}-\ref{algstep-forzi-end} of the $k$-th iteration of the loop for some $k \ge 1$. Assume that the weight
of $z_i$ is set to~$d$ in this iteration.
By Claim (5), we have $g_j \in U_{k}+V_{d-1}$ at this point, so there are $f_j \in U_k$ and
$v_j \in V_{d-1}$ with $g_j = f_j+v_j$.
Hence $f_j = g_j-v_j = z_i-v_j$, where every term in $v_j$ has $W$-weight $\le d-1$, and where $z_i$
has weight~$d$. This shows $\LT_\sigma(f_j) = z_i$ for every term ordering $\sigma$
that is compatible with the grading given by~$W$. Consequently, $(f_1,\dots,f_s)$ is a
$Z$-separating tuple of polynomials.
\end{proof}

Algorithm~\ref{alg-OptimizedCheckZ} sometimes
verifies separating tuples of indeterminates of larger lengths 
than Algorithm~\ref{alg-checkZ}, as the following example shows.

\begin{example}\label{ex-checkZsepI-cont1}
Let $P = \QQ[x_1,\dots,x_{11}]$ and let $I$ be the ideal given in Example \ref{ex-checkZsepI}.
There we used Algorithm \ref{alg-checkZ} to verify that
$(x_4,x_5,x_7)$ is a separating tuple for~$I$.

Now we apply Algorithm~\ref{alg-OptimizedCheckZ} to find a larger separating tuple
of indeterminates for the ideal~$I$. More precisely, consider $Z=(x_4,x_5,x_7,x_9)$.

After deleting all monomials not divisible by indeterminates from~$Z$ in Step~\ref{algstep-delterms2}, we are
left with
\begin{align*}
  g_1 &= x_{4} + x_{1}x_{4}x_{8}x_{11}  + x_{5}x_{6}x_{8} 
        + x_{5}x_{6}x_{10} + x_{7}x_{8}, \quad
  g_2 = x_{7}, \\
  g_3 &= x_{5} + x_{6}x_{7}x_{8} + x_{7}x_{8} + x_{7}x_{10}, \quad
  g_4 = x_{9} -x_{1} x_{7} x_{9}, \\
  g_5 &= x_{1}x_{4}x_{8}x_{11} + x_{5}x_{7}, \quad
  g_6 = x_{7}x_{9}, \quad g_7=g_8=g_9=0.
\end{align*}
At this point, Steps \ref{algstep-interred}a-\ref{algstep-interred}d in Algorithm~\ref{alg-OptimizedCheckZ} 
enable us to reduce $g_4$ against $x_1g_6$ and get the new polynomial
$g_4 = x_9$. Hence we can assign weight 1 to both $x_7$ and $x_9$
in the first iteration of the loop.

The remaining part of the algorithm operates similarly as in Example \ref{ex-checkZsepI}.
However, in this case we find the weight tuple $(0, 0, 0, 157, 13, 0, 1, 0, 1, 0, 0)$
and conclude that the tuple $Z=(x_4,x_5,x_7,x_9)$ is a separating tuple 
of indeterminates for $I$. 

On the other hand, if we apply Algorithm~\ref{alg-checkZ} to $g_1,\dots,g_9$ and $Z$,
then we cannot eliminate the term $x_1 x_7 x_9$ in $\Supp(g_4)$,
which leads to a ``\texttt{Fail}'' in Algorithm~\ref{alg-checkZ}.
\end{example}

Let us finish this discussion with a small optimization of Algorithm~\ref{alg-OptimizedCheckZ}.

\begin{remark} \label{rem-computingBasis-usingDivAlg} 
In the setting of Algorithm~\ref{alg-OptimizedCheckZ}, 
to construct a $K$-basis $\{q_1,\dots,q_u\}$ of the $K$-vector space 
$\langle H \rangle_K \cap P_{\le \delta}$, as required in Steps \ref{algstep-interred}b-\ref{algstep-interred}c, without introducing a degree-compatible
term ordering~$\tau$, one can alternatively use the following instructions:
\begin{itemize}
\item[(1)]
Write $H' = (h_1,\dots,h_m)$ where $H = \{h_1,\dots,h_m\}$ with $m = (n-s)r$. 
Compute the set $\{t_1,\dots,t_e\}$ of all terms of degree 
$> \delta$ in the support of one of the elements in $H$.

\item[(2)]
Write $h_j = h'_{j} +\sum_{k=1}^e a_{jk}t_k$, where 
$a_{jk}\in K$ and $\deg(h'_{j})\le \delta$ for all $j=1,\dots,m$,
and form the matrix $A = (a_{kj}) \in \Mat_{e,m}(K)$.

\item[(3)]
Compute a $K$-basis $\{v_1,\dots,v_{u}\}$ of $\Ker(A)$, and let
$q_k = H\cdot v_k$ for $k=1,\dots,u$.
\end{itemize}

The advantage of this method is that the matrix $A$ may be smaller than the
cofficient matrix of $h_1,\dots,h_m$. This may make the computation of this step more
efficient than the computation in Step~\ref{algstep-interred}b.
\end{remark}

\bigbreak
%
%

\section[Finding Z-Separating Tuples of Polynomials]{Finding $Z$-Separating Tuples of Polynomials}
\label{sec5-compute-polys}

The following algorithm provides a $Z$-separating tuple of polynomials in
the ideal $I$ in the case when Algorithm~\ref{alg-checkZ} is successful.

\begin{algo}{\bf (Computing $Z$-Separating Polynomials)}\label{alg-findSepTuple} \\
Let $g_1,\dots,g_r \in P$ be generators of a proper ideal $I$ in~$P$, and let~$Z$ be a tuple
of distinct indeterminates in~$X$ such that Algorithm~\ref{alg-checkZ}
returns a weight tuple $W \in \NN^n$.
Moreover, let $\sigma$ be a term ordering on~$P$
which is compatible with the grading given by~$W$.
Consider the following sequence of instructions.

\begin{algorithmic}[1]
\State Compute $T=(t_1,\dots,t_m)$ where $\{t_1,\dots,t_m\}$ is the set of terms in the union
of the supports of $g_1,\dots,g_r$, ordered such that $t_1 >_\sigma \dots >_\sigma t_m$. 

\State Form the coefficient matrix $M\in \Mat_{r,m}(K)$ of $g_1,\dots,g_r$ 
with respect to terms in~$T$.

\State Compute the matrix $N \in \Mat_{r',r}(K)$ such that $N \cdot M \in \Mat_{r',m}(K)$ is in
reduced row echelon form and has no zero rows. 

\State Compute $(g_1',\dots,g_{r'}')^{\operatorname{tr}} = N \cdot (g_1,\dots,g_r)^{\operatorname{tr}}$. \;

\State For $i=1,\dots,s$, let $f_i \in \{g_1',\dots,g_{r'}'\}$ be
the element with $\LT_\sigma(f_i) = z_i$. Return $(f_1,\dots,f_s)$.
\end{algorithmic}

\noindent
This is an algorithm which computes a $Z$-separating tuple of polynomials $(f_1,\dots,f_s)$ in~$I$.
\end{algo}

\begin{proof}
First we note that if the algorithm is executed correctly, then the resulting polynomials
$f_1,\dots,f_s$ form indeed a $Z$-separating tuple of polynomials in $I$, 
as they are by construction elements of~$I$ which satisfy $\LT_\sigma(f_i) = z_i$ 
for $i=1,\dots,s$.
It remains to prove that the algorithm is always able to find $f_1,\dots,f_s$
as claimed in Step~4. 

Let $z_i \in Z$. By Proposition \ref{prop-checkZ-success}, 
there is a polynomial $f_i' \in \langle g_1,\dots,g_r \rangle_K$ 
with $\LT_\sigma(f_i') = z_i$.
Since $N$ is an invertible matrix, we also have 
$f_i' \in \langle g_1',\dots,g_{r'}' \rangle_K$,
and since $N \cdot M$ is in reduced row echelon form,
the leading terms of $g_1',\dots,g_{r'}'$ are pairwise distinct.
Hence there exists an index $j \in \{1,\dots,r'\}$ such that 
$\LT_\sigma(g_j') = \LT_\sigma(f_i') = z_i$.
This shows that the algorithm correctly finds $f_i = g_j'$.
\end{proof}

\begin{example}\label{ex-checkZsepI-cont2}
Let $g_1,\dots,g_9$ be the polynomials given in Example~\ref{ex-checkZsepI},
and let $Z = (x_4,x_5,x_7)$. The application of Algorithm~\ref{alg-checkZ}
in Example \ref{ex-checkZsepI} yielded the weight tuple
$W = (0,0,0,43,7,0,1,0,0,0,0)$. As in \cite[Definition~1.4.11]{KR1}, 
we let $\sigma$ be the ordering on $P$ represented by the matrix

\def\mOne{\text{-}1}
$$
A =\left(
\begin{smallmatrix}
  0 & 0 & 0 & 43 & 7 & 0 & 1 & 0 & 0 & 0 & 0 \\
  1 & 1 & 1 & 0 & 0 & 1 & 0 & 1 & 1 & 1 & 1 \\
  0 & 0 & 0 & 0 & 0 & 0 & 0 & 0 & 0 & 0 & \mOne \\
  0 & 0 & 0 & 0 & 0 & 0 & 0 & 0 & 0 & \mOne & 0 \\
  0 & 0 & 0 & 0 & 0 & 0 & 0 & 0 & \mOne & 0 & 0 \\
  0 & 0 & 0 & 0 & 0 & 0 & 0 & \mOne & 0 & 0 & 0 \\
  0 & 0 & 0 & 0 & 0 & 0 & \mOne & 0 & 0 & 0 & 0 \\
  0 & 0 & 0 & 0 & 0 & \mOne & 0 & 0 & 0 & 0 & 0 \\
  0 & 0 & 0 & 0 & \mOne & 0 & 0 & 0 & 0 & 0 & 0 \\
  0 & 0 & \mOne & 0 & 0 & 0 & 0 & 0 & 0 & 0 & 0 \\
  0 & \mOne & 0 & 0 & 0 & 0 & 0 & 0 & 0 & 0 & 0
\end{smallmatrix}
\right).
$$
Then $\sigma$ is a term ordering on P which is compatible with the grading given by~$W$.
An application of Algorithm \ref{alg-findSepTuple} to $g_1,\dots,g_9$, $Z$ and $\sigma$
gives us the $Z$-separating tuple $(f_1,f_2,f_3)$ of polynomials in $I$, where
\begin{align*}
f_1 &= x_{4} -x_{5} x_{7} +x_{5} x_{6} x_{8} +x_{5} x_{6} x_{10} +x_{7} x_{8} -x_{1}^2 x_{6} +x_{3} x_{6} +x_{1}, \\
f_2 &= x_{5} +x_{6} x_{7}^2 x_{8} +x_{7} x_{8} +x_{7} x_{10} +x_{1} x_{2} x_{3}^2 x_{8}^2  +x_{1} x_{3} x_{8}^2  +x_{2}, \\
f_3 &= x_{7} +x_{2}^2 x_{6}^4  +x_{1} x_{2} x_{3}^2 x_{8}^2  +x_{1} x_{2} x_{6}^2 x_{10}^2  +x_{1}^3 x_{10} +x_{3} +1.
\end{align*}
\end{example}

Notice that in the setting of Algorithm~~\ref{alg-OptimizedCheckZ}, it is somewhat trickier to
find actual separating polynomials. Let us describe how to do this.

\begin{remark}\label{rem-findSepTupleOptim}
If Algorithm~\ref{alg-checkZ} is successful, then we can apply
Algorithm~\ref{alg-findSepTuple} in order to find a $Z$-separating tuple of polynomials.
However, if Algorithm~\ref{alg-OptimizedCheckZ} is successful,
then it does not guarantee the existence of a $Z$-separating tuple 
of polynomials in the $K$-vector space $\langle g_1,\dots,g_r \rangle_K$.
Hence we can not apply Algorithm \ref{alg-findSepTuple} to compute such a tuple.
In this case, in order to compute a $Z$-separating tuple of polynomials in~$I$, 
we can mimic all linear operations executed on $g_1,\dots,g_r$ on polynomials
$G_1,\dots,G_r$ as follows. 

At the start of the algorithm, let $G_1,\dots,G_r$ be copies of the input polynomials
$g_1,\dots,g_r$.
In every execution of Steps \ref{algstep-interred}a-\ref{algstep-interred}c, compute the matrix $A \in \Mat_{(n-s)r,u}(K)$ such that
$$
(q_1,\dots,q_u) = (x_ig_j \mid x_i \in X \setminus Z \text{ and } j \in 1,\dots,r) \cdot A
$$
and let
$$
(Q_1,\dots,Q_u) = (x_iG_j \mid x_i \in X \setminus Z \text{ and } j \in 1,\dots,r) \cdot A.
$$
Mimic the linear transformation performed on $g_1,\dots,g_r,q_1,\dots,q_u$ in Step \ref{algstep-interred}d on
$G_1,\dots,G_r,Q_1,\dots,Q_r$, and re-name the resulting polynomials as $G_1,\dots,G_r$ again.
If $x_k = z_i$ in Step \ref{algstep-setweight}, let $f_i = G_k$.

Whenever Algorithm \ref{alg-OptimizedCheckZ} successfully returns a weight tuple~$W$,
we obtain a $Z$-separating tuple $(f_1,\dots,f_s)$ of polynomials in~$I$ in this way.
\end{remark}

The next remark shows how to obtain a coherently $Z$-separating tuple from a
$Z$-separating tuple of polynomials.

\begin{remark}\label{rem-SepToCohSep}
Assume that $(f_1,\dots,f_s)$ is a $Z$-separating tuple of polynomials in~$I$ obtained
by Algorithm \ref{alg-findSepTuple} or Remark \ref{rem-findSepTupleOptim} after a
successful run of Algorithm \ref{alg-checkZ} or \ref{alg-OptimizedCheckZ}, respectively.
Let $W$ be the corresponding output of the latter algorithm
and let $\sigma$ be a term ordering on $P$ which is compatible with the grading given by~$W$.

If we renumber indices such that $z_1 <_\sigma z_2<_\sigma\cdots <_\sigma z_s$,
then $f_i$ is of the form $f_i =z_i-p_i$, where  
$p_i \in K[X\setminus Z][z_1,\dots,z_{i-1}]_{\le \delta}$ for $i=1,\dots,s$.
By substituting $z_1\mapsto p_1$ in $f_2$, we obtain a polynomial
of the from $f'_2=z_2 -p'_{2}$, where $p'_2 \in  K[X\setminus Z]_{\le \delta^2}$.
Next we continue by substituting $z_1 \mapsto p_1$ and $z_2\mapsto p'_2$ in $f_3$ and so on.
In this way, we eventually get polynomials
$f'_1=z_1-h'_1,\dots, f'_s=z_s-h'_s$ with $h'_1,\dots,h'_s\in K[X\setminus Z]$,
i.e., a coherently $Z$-separating tuple $(f'_1,\dots,f'_s)$
of polynomials in $I$ that defines a $Z$-separating re-embedding of~$I$.
\end{remark}

Let us apply this remark to our running example.

\begin{example}\label{ex-checkZsepI-cont3}
Let $P$ and $I$ be defined as in Example~\ref{ex-checkZsepI}, let $Z=(x_4,x_5,x_7,x_9)$,
and let $Y = X \setminus Z$.
In Example \ref{ex-checkZsepI-cont1} we found out that
$Z$ is a separating tuple of indeterminates for~$I$.
Using Remark~\ref{rem-findSepTupleOptim}, we compute
a $Z$-separating tuple $(f_1,\dots,f_4)$ of polynomials in $I$ where
\begin{align*}
f_1 &= x_{7}\; -x_{1}^2 x_{3}^4  -x_{1} x_{2} x_{3} x_{6} x_{8} x_{11} 
+x_{1}^3 x_{10} +x_{3} +1,\\
f_2 &= x_{9}\; +x_{1} x_{8} x_{10}^2  +x_{3} x_{11},\\
f_3 &= x_{5}\; +x_{1} x_{2} x_{3}^2 x_{8}^2  +x_{6} x_{7}^2 x_{8} +x_{1} x_{3} x_{8}^2  
+x_{7} x_{8} +x_{7} x_{10} +x_{2},\\
f_4 &= x_{4} \; - x_{1}^2 x_{6} +x_{5} x_{6} x_{8} +x_{5} x_{6} x_{10} 
+x_{3} x_{6} -x_{5} x_{7} +x_{7} x_{8} +x_{1}.
\end{align*}
Write $p_1 = x_7-f_1$ and $p_2 = x_9-f_2$.
If we substitute $x_7$ with $p_1$ and $x_9$ with $p_2$ in $f_3$, we obtain
$$
f_3' = \; x_5-x_{1} x_{2} x_{3}^2 x_{8}^2  -x_{6}x_{8}p_1^2 
-x_{1} x_{3} x_{8}^2  -x_{8}p_1  -x_{10}p_1 +x_{2}.
$$
Next, substituting $x_7$, $x_9$ and $x_5$ by $p_1$, $p_2$ and $p_3 = x_5-f_3'$,
respectively, into $f_4$ yields
$$
f_4' = x_{4} \, - x_{1}^2 x_{6} +  p_3x_{6} x_{8}
+ p_3x_{6}x_{10} +x_{3} x_{6} - p_3p_1 + p_1x_{8} +x_{1}. 
$$
Altogether, the tuple $(f_1,f_2,f_3',f_4')$ is a coherently $Z$-separating tuple
of polynomials in~$I$ which defines a $Z$-separating re-embedding of $P/I$ into
$K[Y]/(I\cap K[Y])$.
If we now substitute $x_4$, $x_5$, $x_7$ and $x_9$ with $p_1$, $p_2$, $p_3$, and
$p_4 = x_4-f_4'$, respectively, in the generators $g_1,\dots,g_9$ of $I$,
we obtain polynomials $\bar{g}_1,\dots,\bar{g}_9 \in K[Y]$ which are a generating
set of $I \cap K[Y]$.
\end{example}

\bigbreak
%
%

\section{Checking Separating Boolean Indeterminates}
\label{sec6-boolpolys}

Throughout this section we let $\FF_2 = \ZZ/2\ZZ$, let
$P = \FF_2[X_1,\dots,X_n]$ be the polynomial ring over $\FF_2$ in the indeterminates
$X_1,\dots,X_n$, and let $\BB_n = P/\II_n$ be the \textbf{ring of Boolean polynomials}. 
Here $\II_n = \langle X_1^2-X_1,\dots,X_n^2-X_n\rangle$ is the \textbf{field ideal} of~$P$. 
For $i \in \{1,\dots,n\}$, we denote the residue class of $X_i$ in $\BB_n$ by $x_i$.
These residue classes $x_1,\dots,x_n$ are also referred to as \textbf{Boolean indeterminates}.
\medskip

The goal of this section is to introduce (coherently) separating Boolean polynomials
and to devise algorithms for checking whether given tuples of Boolean indeterminates
are separating, and in the positive case for computing the corresponding separating Boolean polynomials.
For basic facts about Boolean polynomials, 
we refer the reader to~\cite{Bri} or~\cite[Sec.~2.10]{Hor}.
\medskip

It is well-known (see e.g.~\cite[Theorem 2.1.2]{Bri}) that every Boolean polynomial
$f \in \BB_n$ has a unique representative $F \in P$ that is a sum of square-free terms,
i.e., terms of the form $X_1^{\alpha_1}\cdots X_n^{\alpha_n}$ with $\alpha_i\in \{0,1\}$.
In this case, the polynomial $F$ is called the \textbf{canonical representative}
of~$f$. We write $\deg(f) = \deg(F)$ and $\Supp(f) = \{ t+\II_n \mid t \in \Supp(F) \}$.
Moreover, given a term ordering~$\sigma$ on~$P$, we let
$\LT_\sigma(f) = \LT_\sigma(F)+\II_n$
and call it the {\bf leading term} of the Boolean polynomial~$f$ with respect to~$\sigma$.
We also write $\NR_{\sigma,G}(f) = \NR_{\sigma,(G_1,\dots,G_r,X_1^2-X_1,\dots,X_n^2-X_n)}(F)+\II_n$ 
for the normal remainder of~$f$
with respect to a tuple $G = (g_1,\dots,g_r) \in \BB_n^r$, where $G_i$ is the canonical
representative of $g_i$ for $i=1,\dots,r$.
This allows us also to define a \textbf{Boolean $\sigma$-Gr\"obner basis} of~$I$ 
as a finite set $G \subseteq I$ that satisfies
$\langle \LT_\sigma(g) \mid g \in G \rangle = \LT_\sigma(I) 
= \langle \LT_\sigma(f) \mid f \in I \rangle$.
\medskip

In the following, let $X = (x_1,\dots,x_n)$ be the tuple of all 
Boolean indeterminates in~$\BB_n$, and let $Z = (z_1,\dots,z_s)$ be a tuple 
of distinct Boolean indeterminates from~$X$.
In analogy to Definition~\ref{def-Zsep}, the tuple $Z$ is called a \textbf{separating tuple
of Boolean indeterminates} for a proper ideal~$I$ in~$\BB_n$, if there exist a term
ordering $\sigma$ and a tuple $F = (f_1,\dots,f_s)$ of Boolean polynomials
in~$I$ such that $\LT_\sigma(f_i) = z_i$ for $i=1,\dots,s$.
The tuple $F$ is then called a \textbf{$Z$-separating tuple of Boolean polynomials} in~$I$.
Whenever such a tuple~$F$ satisfies the additional condition that none of the terms in the support 
of any~$f_i$ is divisible by $z_j$ for $i,j \in \{1,\dots,s\}$ with $i \neq j$,
we say that~$F$ is a \textbf{coherently $Z$-separating tuple of Boolean polynomials} in~$I$.
\medskip

The connection between $Z$-separating and coherently $Z$-separating tuples of
Boolean polynomials is described by the following proposition.
The analogous proposition for usual polynomials can be found as Proposition~2.2.d
in~\cite{KLR2}.

\begin{proposition}\label{prop-cohSepPoly-bool-1}
Let $I$ be a proper ideal in $\BB_n$, let $(f_1,\dots,f_s)$ be a $Z$-separating
tuple of Boolean polynomials in $I$, and let $\sigma$ be a term ordering on $\BB_n$ with
$\LT_\sigma(f_i) = z_i$ for $i=1,\dots,s$.
Moreover, write $f_i = z_i+h_i$ with $h_i \in \BB_n$, and let
$\tilde{h}_i = \NR_{\sigma, (f_1,\dots,f_s)}(h_i)$ for $i=1,\dots,s$.
Then the tuple $(z_1-\tilde{h}_1,\dots,z_s-\tilde{h}_s)$
is a coherently $Z$-separating tuple of Boolean polynomials in $I$.
\end{proposition}

\begin{proof}
For $i \in \{1,\dots,s\}$, let $F_i, Z_i, \widetilde{H}_i$ be the
canonical representatives of $f_i,z_i,\tilde{h}_i$, respectively.
By the definition of the normal remainder, the polynomials $\widetilde{H}_i$ 
are fully reduced against $F_1,\dots,F_s$. Thus no term in their support
is divisible by any of the leading terms of $F_1,\dots,F_s$. 
Since $\LT_\sigma(F_j)=Z_j$ for $j=1,\dots,s$, this implies that $Z_j$ does not divide 
any of the terms in the support of $\widetilde{H}_i$,
and hence $z_j$ does not divide any of the terms in the support of $\tilde{h}_i$
for $i,j\in \{1,\dots,s\}$. 
Consequently, none of the terms in the support of $z_i-\tilde{h}_i$ is divisible by $z_j$
for $i,j \in \{1,\dots,s\}$ with $i \neq j$. 
In addition, it is easy to see that $h_i-\tilde{h}_i 
\in \langle f_1,\dots,f_s \rangle \subseteq I$, and so
$z_i-\tilde{h}_i = f_i- (h_i-\tilde{h}_i) \in I$ for $i=1,\dots,s$.
Therefore the tuple $(z_1-\tilde{h}_1,\dots,z_s-\tilde{h}_s)$ is a coherently
$Z$-separating tuple of Boolean polynomials in $I$.
\end{proof}

In order to construct separating re-embeddings of an ideal $I$ in $\BB_n$,
one can apply the result below.
Here a {\bf Boolean re-embedding} of $I$ is an $\FF_2$-algebra isomorphism 
$$
\Psi:\; \BB_n/I \;\longrightarrow\; \BB_m/I'
$$
where $m \in \NN$ and $I'$ is an ideal of~$\BB_m$.

\begin{proposition}\label{prop-cohSepPoly-bool}
{\bf (Boolean $Z$-Separating Re-Embeddings)} \\
In the setting described above, let~$I$ be a proper ideal in~$\BB_n$, and let $\BB_m = \FF_2[X\setminus Z]$
with $m=n-s$.
Assume that $F = (f_1,\dots,f_s)$ is a coherently $Z$-separating tuple 
of Boolean polynomials in~$I$ with respect to a term ordering~$\sigma$ on~$\BB_n$.
For $i=1,\dots,s$, write $f_i = z_i -h_i$ with $h_i \in \BB_m$.  
\begin{enumerate}
\item[(a)] The reduced Boolean $\sigma$-Gr\"obner basis of $I$ is of the form 
$$
G = \{z_1-\tilde{h}_1,\dots,z_s-\tilde{h}_s,g_1,\dots, g_r\},
$$
where $\tilde{h}_1,\dots,\tilde{h}_s, g_1,\dots,g_r \in \BB_m$. 

\item[(b)] The $\FF_2$-algebra homomorphism $\Phi: \BB_n/I \longrightarrow \BB_m / 
(I ~\cap~ \BB_m)$ given by 
$$
\Phi(x_i+I) \;=\; \begin{cases}  x_i+(I \cap \BB_m) & \text{ for }x_i\notin Z\\  
                            h_j+(I \cap \BB_m) & \text{ for }x_i = z_j \in Z
                  \end{cases}
$$ 
is an isomorphism of $\FF_2$-algebras.
It is called a \textbf{Boolean $Z$-se\-pa\-ra\-ting re-embedding} of~$I$.

\item[(c)] For every elimination ordering $\tau$ for~$Z$, we have 
$\langle Z\rangle \subseteq \LT_\tau(I)$.
\end{enumerate}
\end{proposition}

\begin{proof}
Let $G$ be the reduced Boolean $\sigma$-Gr\"obner basis of~$I$, and let 
$\tilde{h}_i = \NR_{\sigma, G}(h_i)$ for $i=1,\dots,s$.
By assumption, we have $z_i = \LT_\sigma(f_i) \in \LT_\sigma(I)$ 
and $h_i\in \BB_m$. In particular, since $\LT_\sigma(I)\ne \langle 1 \rangle$, we have 
$z_i-\tilde{h}_i\in G$ and $\tilde{h}_i\in \BB_m$ for $i=1,\dots,s$. 
Thus we can write $G = \{z_1-\tilde{h}_1,\dots,z_s-\tilde{h}_s,g_1,\dots, g_r\}$,
where $g_1,\dots,g_r \in \BB_n$. Moreover, since~$G$ is a reduced Gr\"obner basis, we know 
that~$z_i$ does not divide any term in $\Supp(g_j)$ for $i\in \{1,\dots,s\}$ and $j\in\{1,\dots,r\}$.
Hence we get $g_1,\dots,g_r \in \BB_m$, and claim~(a) follows.

To show~(b), it suffices to prove that $I=\langle G \rangle$ is the kernel of the 
$\FF_2$-algebra epimorphism $\phi: \BB_n \longrightarrow \BB_m/(I \cap \BB_m)$ defined by 
$x_i \mapsto x_i+(I \cap \BB_m)$ for $x_i \notin Z$ and by $x_i \mapsto h_j+(I \cap \BB_m)$ 
for $x_i = z_j \in Z$.

Clearly, $G\subseteq \Ker(\phi)$, and so $I=\langle G \rangle\subseteq \Ker(\phi)$.
For the other inclusion, let $h \in \Ker(\phi)$ and write $\tilde{h} = \NR_{\sigma,G}(h)$.
Then $\tilde{h} \in \BB_m$ and $h = f+\tilde{h}$ for a Boolean polynomial $f \in I$.
Since we already know $I \subseteq \Ker(\phi)$, this implies 
$\phi(\tilde{h}) = \phi(h) = 0$, and so $\tilde{h} \in I\cap \BB_m$.
Hence we get $h = f+\tilde{h} \in I$.

Finally, claim (c) follows from the observation that 
$\LT_\tau(f_i)=\LT_\tau(z_i-h_i) = z_i \in \LT_\tau(I)$ for
$i=1,\dots,s$, because~$\tau$ is an elimination ordering for~$Z$.
\end{proof}

The next proposition provides a connection between $Z$-separating Boolean polynomials 
and $Z$-separating polynomials in~$P$.

\begin{proposition}\label{prop-boolZsepConnection}
Let $I'$ be a proper ideal in $P$ that contains the field ideal $\II_n$, and let $I = I'/\II_n$.
Let $Z = (z_1,\dots,z_s)$ be a tuple of distinct Boolean indeterminates, and let
$Z' = (Z_1,\dots,Z_s)$ be the tuple of their canonical representatives, i.e., let
$Z_i \in \{X_1,\dots,X_n\}$ and $z_i = Z_i+\II_n$ for $i=1,\dots,s$.
Furthermore, let $F_1,\dots,F_s \in P$ and let $f_i = F_i+\II_n$ for $i=1,\dots,s$.
\begin{enumerate}
\item[(a)] If $(F_1,\dots,F_s)$ is a $Z'$-separating tuple of polynomials in $I'$,
then $(f_1,\dots,f_s)$ is a $Z$-separating tuple of Boolean polynomials in $I$.

\item[(b)] If $(f_1,\dots,f_s)$ is a $Z$-separating tuple in $I$ and $F_1,\dots,F_s$
are the canonical representatives of $f_1,\dots,f_s$, respectively,
then $(F_1,\dots,F_s)$ is a $Z'$-separating tuple of polynomials in $I'$.

\item[(c)] The tuple $Z'$ is a separating tuple of indeterminates for $I'$ if and only if 
the tuple~$Z$ is a separating tuple of Boolean indeterminates for~$I$.

\end{enumerate}
\end{proposition}

\begin{proof}
To show (a), let $\sigma$ be a term ordering on $P$ such that $\LT_\sigma(F_i) = Z_i$ for
$i=1,\dots,s$.
Fix $i \in \{1,\dots,s\}$. Consider $t = x_{i_1} \cdots x_{i_k} \in \Supp(f_i)\setminus\{z_i\}$, 
where $k=\deg(t)$. Then there are $\alpha_1,\dots,\alpha_k \in \NN_+$ with
$T = X_{i_1}^{\alpha_1} \cdots X_{i_k}^{\alpha_k} \in \Supp(F_i)$.
It follows from $\LT_\sigma(F_i) = Z_i$ that 
$Z_i = \LT_\sigma(F_i) >_\sigma T \geq_\sigma X_{i_1} \cdots X_{i_k}$,
and hence $z_i >_\sigma x_{i_1} \cdots x_{i_k} = t$.
This implies that $\LT_\sigma(f_i) = z_i$. Therefore $(f_1,\dots,f_s)$ is a
$Z$-separating tuple of Boolean polynomials in~$I$.

To prove (b), let $\sigma$ be a term
ordering such that $\LT_\sigma(f_i) = z_i$ for $i=1,\dots,s$.
Then by the definition of leading terms of Boolean polynomials, we have
$\LT_\sigma(F_i) = Z_i$ for all $i=1,\dots,s$, and thus the tuple $(F_1,\dots,F_s)$
is a $Z'$-separating tuple of polynomials in $I'$.

By the definition of separating tuples of indeterminates, (c) follows immediately from~(a) 
and~(b).
\end{proof}

After these preparations, we are now ready to adapt our algorithms to 
check whether a tuple $Z$ of Boolean indeterminates 
is a separating tuple for a given ideal in~$\BB_n$. 
The following result plays an important role for that.

\begin{proposition}[Computation of $Z$-Separating Boolean Polynomials]\label{prop-checkZ-bool}
$\mathstrut$\\
Let $I$ be a proper ideal in $\BB_n$, let $\{g_1,\dots,g_r\} \subseteq \BB_n$ 
be a system of generators of~$I$, and let $G_1,\dots,G_r\in P$ be the
canonical representatives of $g_1,\dots,g_r$.
Furthermore, let $Z = (z_1,\dots,z_s)$ be a tuple of distinct Boolean indeterminates, 
and let $Z'=(Z_1,\dots,Z_s)$ be the tuple of their canonical representatives.
\begin{itemize}
\item[(a)] If we apply Algorithm \ref{alg-checkZ} to the input $(G_1,\dots,G_r)$ and~$Z'$,
it successfully returns a weight tuple $W \in \NN^n$ if and only if
there exist $f_1,\dots,f_s \in \langle g_1,\dots,g_r \rangle_{\FF_2}$
such that $(f_1,\dots,f_s)$ is a $Z$-separating tuple of
Boolean polynomials in~$I$.

\item[(b)] In the case of~(a), let $\sigma$ be any term ordering which is compatible with
the grading given by~$W$. Then there exist
$f_1',\dots,f_s' \in \langle g_1,\dots,g_r \rangle_{\FF_2}$
with $\LT_\sigma(f_i') = z_i$ for $i=1,\dots,s$.

\end{itemize}
\end{proposition}

\begin{proof}
Let $U = \langle g_1,\dots,g_r \rangle_{\FF_2}$ and
let $U' = \langle G_1,\dots,G_r \rangle_{\FF_2}$.
Consider the $\FF_2$-linear map  $\phi : \BB_n \to P$ that maps every 
Boolean polynomial to its canonical representative. 
Then the restriction map $\phi|_U: U\rightarrow U'$ 
is an isomorphism of $\FF_2$-vector spaces.
We first show that there is a $Z$-separating tuple of Boolean polynomials for $I$
in $U$ if and only if there is a~$Z'$-separating tuple of polynomials
for $I' = \langle G_1,\dots,G_s \rangle + \II_n$ in $U'$.

Suppose that $(f_1,\dots,f_s) \in U^s$ is a $Z$-separating tuple of Boolean polynomials for~$I$.
Then $(\phi(f_1),\dots,\phi(f_s)) \in (U')^s$ is a $Z'$-separating tuple for~$I'$
by Proposition~\ref{prop-boolZsepConnection}.b.
Conversely, let $(F_1,\dots,F_s) \in (U')^s$ be a $Z'$-separating
tuple of polynomials for~$I'$, and let $f_i = F_i+\II_n$ for $i=1,\dots,s$.
Then $(f_1,\dots,f_s)$ is a $Z$-separating tuple for $I$ by
Proposition~\ref{prop-boolZsepConnection}.a.
Since $F_1,\dots,F_s \in U'$, every term in their support is 
in the support of one of the polynomials $G_1,\dots,G_r$ and thus square-free. 
It follows that every $F_i$ is the canonical representative of $f_i$
and therefore $f_i = \phi^{-1}(F_i)\in \phi^{-1}(U') = U$ for $i=1,\dots,s$.

In Proposition~\ref{prop-checkZ-success} it was shown that Algorithm~\ref{alg-checkZ}
successfully returns a weight tuple~$W$ if and only if there exists a
$Z'$-separating tuple of polynomials in~$U'$.
By the observation above, this is equivalent to the condition that there exists 
a $Z$-separating tuple of polynomials for $I$ in $U$.
This shows (a).

For (b), let $\sigma$ be a term ordering which is compatible
with the grading given by $W$. By Algorithm~\ref{alg-checkZ}
and Proposition~\ref{prop-checkZ-success},
there exists a $Z'$-separating tuple $(F_1',\dots,F_s')$ in~$U'$ 
such that $\LT_\sigma(F'_i)=Z_i$ for $i=1,\dots,s$. 
The canonical representative of the element $f_i' = F_i'+\II_n$ 
is $F'_i$ for $i=1,\dots,s$. Therefore we conclude that
$\LT_\sigma(f_i') = \LT_\sigma(F_i')+\II_n = z_i$ 
for $i=1,\dots,s$.
\end{proof}

\begin{remark}\label{rem-checkZ-bool}
Let $Z = (z_1,\dots,z_s)$ be a tuple of distinct Boolean indeterminates in~$X$,
let $g_1,\dots,g_r \in \BB_n$, and let $Z_i,G_j$ be the canonical
representatives of $z_i,g_j$, respectively, where $i\in \{1,\dots,s\}$ 
and $j\in \{1,\dots,r\}$.
The following methods can be applied to check whether $Z$ is a separating tuple 
for the ideal $I = \langle g_1,\dots,g_r \rangle$.
\begin{itemize}
\item[(a)]
Apply Algorithm \ref{alg-checkZ} to $G_1,\dots,G_r$. 
By Proposition \ref{prop-checkZ-bool}.a,
Algorithm \ref{alg-checkZ} successfully returns a weight tuple if and only
if $Z$ is a separating tuple of Boolean indeterminates for~$I$ 
and there is a $Z$-separating tuple of polynomials for $I$ in 
$\langle g_1,\dots,g_r \rangle_K$.

\item[(b)] 
Suppose that Algorithm \ref{alg-checkZ} successfully returns 
a weight tuple $W \in \NN^n$, and let $\sigma$ be a term ordering 
which is compatible with the grading given by $W$.
Apply Algorithm~\ref{alg-findSepTuple} to $(G_1,\dots,G_r)$
and $(Z_1,\dots,Z_s)$ to obtain polynomials $F_1,\dots,F_s \in P$.
Then we have $F_1,\dots,F_s \in \langle G_1,\dots,G_r \rangle_{\FF_2}$
with $\LT_\sigma(F_i) = Z_i$ for $i=1,\dots,s$. Since all terms in the
supports of $G_1,\dots,G_r$ are square-free, all terms in the supports
of $F_1,\dots,F_s$ are also square-free. 
Thus $F_i$ is the canonical representative of $f_i = F_i+\II_n$, 
and hence $\LT_\sigma(f_i) = \LT_\sigma(F_i)+\II_n = z_i$.
This implies that $(f_1,\dots,f_s)$ is a $Z$-separating tuple of 
Boolean polynomials for $I$.

\item[(c)]
Let $\{H_1,\dots,H_n\}$ be a set of generators of~$\II_n$,
let $I' = \langle G_1,\dots,G_r \rangle + \II_n$, and let $Z'=(Z_1,\dots,Z_s)$.
According to Proposition \ref{prop-boolZsepConnection}.c, if 
Algorithm~\ref{alg-OptimizedCheckZ} applied to $G_1,\dots,G_r$,
$H_1,\dots,H_n$ and $Z'$ verifies that $Z'$ is a separating tuple 
of indeterminates for $I'$, then the tuple~$Z$ is a separating tuple 
of Boolean indeterminates for~$I$.
In this case, a $Z'$-separating tuple $(F_1,\dots,F_s)$ of polynomials in~$I'$ 
can be computed by using Remark \ref{rem-findSepTupleOptim}. Then 
Proposition \ref{prop-boolZsepConnection}.a yields that 
$(f_1,\dots,f_s)$, where $f_i = F_i+\II_n$ for $i=1,\dots,s$,
is a $Z$-separating tuple of polynomials in~$I$.

\item[(d)]
Suppose that there exist $j \in \{1,\dots,r\}$ and $i \in \{1,\dots,s\}$ 
such that $g_j\in I$ has no constant term and $z_i\in \Supp(g_j)$. 
Then the linear part of $z_ig_j$ is $z_i$. In some cases, it may be useful 
to add those polynomials $z_ig_j$ to the generators of the ideal to 
find separating tuples of Boolean indeterminates.

\end{itemize}
\end{remark}

\bigbreak
%
%

\section{Application to the Cryptoanalysis of AES-128}
\label{sec7-appl}

The Advanced Encryption Standard (AES) was certified by the
U.S.~National Institute of Standards and Technology (NIST) in 2001 
(see~\cite{aes}) and is one of the most widely used symmetric encryption ciphers.
In this section we apply our algorithms to polynomial
equations describing AES-128 with the goal of improving algebraic attacks.
We assume that the reader is familiar with the basic definition of~AES, 
for instance as laid out in~\cite{aes}.

\begin{example}[Polynomial Representation of the AES $S$-Box]\label{ex-sbox-polys}$\mathstrut$\\
The only nonlinear part during the encryption process of AES-128 is the application
of the Rijndael $S$-Box. This is
a map $s : \FF_2^8 \to \FF_2^8$ which is detailed in Section 5.1.1 of \cite{aes}.
We describe~$s$ by the vanishing ideal $I_S$ of its graph, i.e., by the ideal~$I_S$ in~$\BB_{16}$ 
in the Boolean indeterminates $x_1,\dots,x_8,y_1,\dots,y_8$ which is the vanishing ideal of the set
$S = \{(a,s(a)) \mid a \in \FF_2^8\} \subseteq \FF_2^{16}$.

\begin{itemize}
\item[(a)]
Using a computer algebra system such as \cocoa-5 (see \cite{CoCoA}), we can start from the
256 points of~$S$ and use them to compute a generating set of~$I_S$.
In our case we used the function \texttt{IdealOfPoints} of \cocoa-5, which resulted in a
set of~$39$ quadratic Boolean polynomials. We reference this set
as $\Gstandard$ for the remaining part of this section.

\item[(b)]
Now we apply our implementation of Algorithm \ref{alg-checkZ} in the computer algebra system
\texttt{SageMath} (see \cite{sage}) to the polynomials in $\Gstandard$ and every
$Z \subseteq \{x_1,\dots,x_8\}$.
The result is that Algorithm \ref{alg-checkZ} is successful for every~$Z$ with $\#Z \le 3$.
We choose $Z' = (x_1,x_2,x_3)$ and use Algorithm \ref{alg-findSepTuple} to
compute a $Z'$-separating tuple $G_{123}$. The polynomials in $G_{123}$ are of degree~2.

\item[(c)]
Proceeding as in~(b), we also apply our implementation of Algorithm \ref{alg-OptimizedCheckZ} to
$\Gstandard$ and every $Z \subseteq \{x_1,\dots,x_8\}$. It turns out to be successful for every
tuple~$Z$ with $\#Z \le 6$. We choose $Z'' = (x_1,x_2,x_3,x_4,x_5,x_6)$ 
and use Remark \ref{rem-findSepTupleOptim}
to compute a $Z''$-separating tuple $\G6$ in~$I_S$. The polynomials in~$\G6$ are of degree~3.

\item[(d)] If we strive to find a separating system of generators of~$I_S$ consisting of polynomials
of lowest possible degree, we can take the three polynomials of degree~2 in~$G_{123}$, 
three polynomials of degree~3
in~$\G6$, and two polynomials of degrees~5 and~7, respectively, found using further rounds 
of~Algorithm~\ref{alg-OptimizedCheckZ} or suitable elimination computations.
The resulting tuple $G_{1-8}$ is separating for $(x_1,\dots,x_8)$ and generates~$I_S$.

\item[(e)] Finally we can compute a ${\tt lex}$-Gr\"obner basis of~$I_S$ and get a
$(x_1,\dots,x_8)$-separating tuple $G_{lex}$ of polynomials of degree~7.
\end{itemize}

The sets $\Gstandard$, $\Gstandard\cup G_{123}$, $\Gstandard\cup\G6$, $G_{1-8}$ and $G_{lex}$ 
generate $I_S$ and allow us to eliminate 0, 3, 6, 8, and 8 Boolean indeterminates, 
respectively, from~$I_S$. To perform the experiments below, we have to convert these
Boolean polynomials to a suitable input for the solvers we consider. While \bosphorus{}
and {\tt PolyBoRi} work directly with algebraic normal forms and do not require any conversions,
the solver \kissat{} needs CNF (conjunctive normal form) input, the solver \cms{} expects
CNF-XOR, and the solver {\tt 2-Xornado} uses XNF (xor-or-and normal form) formulas.
Let us briefly collect the numbers of logical variables and clauses produced by the 
conversion tools described in~\cite{ADK}.

\begin{table}[h!]
\centering
\caption{Logic Conversions}
\begin{tabular}{cccc}
  \hline
  \multirow{2}{*}{Generating set of $I_S$} & \multicolumn{3}{c}{\# variables/clauses} \\
  & CNF & CNF-XOR & XNF \\
  \hline
  $\Gstandard$ &  \,780/8722 & 256/519  & 136/399  \\
  $G_{lex}$    &  1062/6488  & 843/1424 & 555/1136 \\
  $G_{1-8}$   &  \,542/4061 & 346/606  & 213/473 \\
  \hline
\end{tabular}
\footnotetext{Numbers of logical variables and clauses for the conversions of 
different systems of generators of~$I_S$.}
\label{tab:conversion}
\end{table}

Of course, it does not make sense to include the numbers for $\Gstandard\cup G_{123}$
and $\Gstandard \cup \G6$ since these systems of generators properly contain~$\Gstandard$.
In all three cases, the conversion of $G_{lex}$ is the least efficient one, while
most of the numbers for $\Gstandard$ and $G_{1-8}$ are similar. However, the number of CNF clauses
for $\Gstandard$ is large, because the conversion process creates long linear forms which
blow up the CNF. An adjusted splitting technique (see~\cite{JK}) could possibly improve this
to some extent.
\end{example}

The question studied below is whether adding the polynomials
in $G_{123}$ and $\G6$ speeds up algebraic attacks on AES-128.
Let us describe the mathematical setting of general algebraic attacks first.

\begin{remark}[Algebraic Attacks on $r$-AES-128]\label{rem-AES-polys}$\mathstrut$\\
The original cipher AES-128 is specified to have 10 rounds.
For $r \in \{1,\dots,10\}$, we denote by $r$-AES-128 the version of
AES-128 with $r$ instead of $10$ rounds, i.e., the number \texttt{Nr} in Figure 5
of \cite{aes} is set to~$r$.

Consider the ring of Boolean polynomials $\BB_n$ where $n = 128 \cdot (4r+3)$ whose
indeterminates are $\k_i$, $\p_i$, $\c_i$, $\rk_{j,i}$, $\si_{j,i}$, $\so_{j,i}$, and $\rks_{j,i}$
for $i=1,\dots,128$ and $j=1,\dots,r$.
Here the indeterminates $\k_i$, $\p_i$, $\c_i$, and $\rk_{j,i}$
represent key, plaintext, ciphertext, and round key bits, respectively.
The indeterminates $\si_{j,i}$ and $\so_{j,i}$ are additional variables for the inputs and outputs of the
\texttt{Subbytes} transformations (see Section~5.1.1 in \cite{aes})
during the encryption process (see Figure 5 in \cite{aes}).
Lastly, the inderterminates $\rks_{j,i}$ are additional variables for the outputs of the \texttt{Subbytes}
transformations during the computation of the round keys (see Figure 11 in \cite{aes}).

Let $I_r \subseteq \BB_n$ be the ideal modelling $r$-AES-128. In other words, a tuple
$(\k,\p,\c,\rk,\si,\so,\rks) \in (\FF_2^{128})^3 \times (\FF_2^{128r})^4$
is a zero of~$I_r$ if and only if $\p$ encrypted with $r$-AES-128 using the key $\k$
yields~$\c$, and during the encryption process, the round keys are $\rk$, the
$S$-box inputs and outputs are $\si$ and~$\so$, respectively, and the
outputs of the $S$-boxes during the round key computations are given by~$\rks$.
Using a generating set $G$ for $I_S$ as in Example~\ref{ex-sbox-polys},
it is now possible to construct a generating set~$H$ of $I_r$ that only contains
linear polynomials and polynomials from~$G$ with their indeterminates substituted by
the corresponding input and output indeterminates of the \texttt{Subbytes} steps.
\end{remark}

The purpose of the next example is to verify that known-plaintext attacks on AES-128
based on finding the zeros of a specialization of the ideal $I_r$ constructed
in the preceding remark become easier if we enlarge the sets of generators by $G_{123}$
or $\G6$ to eliminate some indeterminates in the $S$-boxes.
Currently, the most efficient tools to solve such systems are SAT solvers.
In \cite{ADK} it is described how systems of Boolean polynomials 
can be converted to XNF (xor-or-and normal form), to XOR-CNF, and to CNF
such that they become suitable inputs for current SAT solvers.

\begin{example}[Attacking One Round of AES-128]\label{ex-1-AES} $\mathstrut$\\
For the following experiments, we generated Boolean polynomial systems~$H$ representing $1$-AES-128.
In order to use them to perform a known-plaintext attack, we proceeded as follows.
\begin{enumerate}
\item[(a)] First we constructed a plaintext-ciphertext pair $(p,c) \in\FF_2^{256}$
for 1-AES-128, where $p = (p_1,\dots,p_{128})$ and $c = (c_1,\dots,c_{128})$.

\item[(b)] Then we specialized various systems of generators~$H$ of the ideal~$I_1$
defining the cipher by letting $\p_i \mapsto p_i$ and $\c_i \mapsto c_i$
for $i=1,\dots,128$. Here the various systems~$H$ differ by the system of generators~$G$
which is used to represent the $S$-box.

\item[(c)] Next we transformed the specializations to suitable inputs
for state-of-the-art Boolean system solvers.

\item[(d)] Finally, we used these solvers to determine the unknown key bits $\k_i$
in 100 randomly generated instances and measured their average solving time.
\end{enumerate}

In Table \ref{tab:benchmarks}, we compare the running times of the XOR-CNF SAT solver
\cms{}, the CNF SAT solver \kissat{} (see \cite{BD}), and the
Boolean polynomial solver \bosphorus{} (see \cite{CSCM}).
For \cms{} and \kissat{}, we converted the Boolean polynomial
systems~$H$ to XOR-CNF and CNF, respectively, using the methods described in~\cite{ADK}.

\begin{table}[h!]
\centering
\caption{Experiments}
\begin{tabular}{cccc}
  \hline
  \multirow{2}{*}{Generating set of $I_S$} & \multicolumn{3}{c}{Average solving time} \\
  & \cms{} & \kissat{} & \bosphorus{} \\
  \hline
  $\Gstandard$               & 40.23 s & 50.53 s & 244.74 s \\
  $G_{123} \cup \Gstandard$  & 35.69 s & 27.36 s & 257.01 s \\
  $\G6 \cup \Gstandard$      & 21.47 s & 33.57 s & 252.78 s \\
  \hline
\end{tabular}
\footnotetext{Average solving times for 100 polynomial systems coming from
  known-plaintext attacks on one round of AES-128.}
\label{tab:benchmarks}
\end{table}

As can be seen from this table, the Boolean polynomial system solver \bosphorus{} is 
significantly slower than the SAT solvers and is not effected much by the choice of the
generating set of~$I_S$. Both \cms{} and \kissat{} profit substantially from the
addition of the separating tuples of polynomials $G_{123}$ and $\G6$.
More precisely, \cms{} performs best when using $\G6 \cup \Gstandard$ and
\kissat{} performs best when using $G_{123} \cup \Gstandard$ instead of $\Gstandard$.

These Boolean polynomial systems were also given to the algebraic solver {\tt PolyBoRi} (see~\cite{BD})
and to the XNF-solver \texttt{2-Xornado} (see~\cite{ADK}). However, neither of these solvers 
terminated on any of the instances in under one hour.
\end{example}

For our next experiment, we considered 2-AES-128 and assumed that a certain 
amount of bits of the secret key was known, as happens for instance in the case of
some sidechannel attacks or for some types of guess-and-determine attacks.

\begin{example}[Attacking Two Rounds of AES-128]\label{ex-2-AES}$\mathstrut$\\
For this experiment, we created 220 instance of algebraic attacks to 2-AES-128.
As above, we substituted randomly generated known plaintext-ciphertext pairs
into the Boolean polynomials. Then we constructed 20 instances for each
value $k\in \{70,\dots,80\}$, where we assume that we know the first~$k$ bits 
of the secret key. These values were also chosen randomly and substituted
into the equations.

In the following cactus plot, we draw the solving times that \bosphorus{} and \cms{}
needed for several systems of generators~$H$ of 2-AES-128 derived from the
stated systems of generators~$G$ of the vanishing ideal~$I_S$ of the S-boxes.

\begin{figure}[H]
  \centering
  \includegraphics[width=0.85\textwidth]{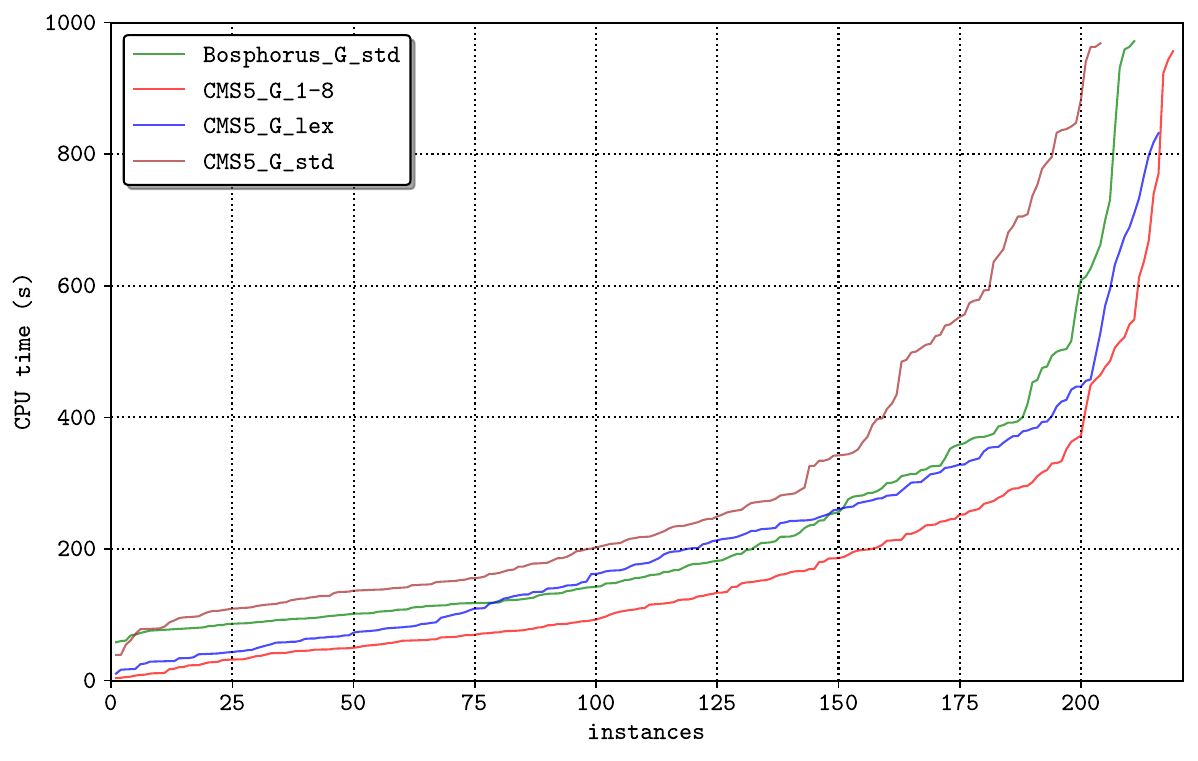}
  \caption{
    Running times of \{{\tt solver}\}\_\{{\tt gens}\} for 
    220 satisfiable instances related to key-recovery attacks on 2-AES-128 
    with knowledge of the first $70\le k\le 80$ key bits.
  }
  \label{fig:benchmarks-2rounds}
\end{figure}

As one can see, the Boolean polynomial system based on the polynomials $G_{1-8}$ constructed
using Algorithms~\ref{alg-checkZ} and~\ref{alg-OptimizedCheckZ} is solved significantly
faster by \cms{} (denoted by {\tt CMS5}) than the system based on the original generators~$\Gstandard$. 
The solver \bosphorus{} did not terminate on the instances involving $G_{lex}$ and~$G_{1-8}$. 
\end{example}

All experiments in this section were performed on an AMD Ryzen 5 7530U CPU 
with 16 GB of RAM under Manjaro Linux.
We used \cms{} version 5.11.22, \kissat{} version 4.0.1, and
\bosphorus{} version 3.0.0.

%
%

\section*{Conclusion}

In recent years, a new technique for performing elimination on large-scale
examples of polynomial ideals has emerged. It has been called elimination by substitution
or the method of $Z$-separating re-embeddings. In this paper we made another
important step in the development of this technique. Previous works discussed 
fast ways of finding tuples of indeterminates~$Z$ which are candidates for this
type of elimination computation~\cite{KLR2,KLR3,KR3}. Frequently there is a large 
number of such candidate tuples. Here we provide several
algorithms which allow us to check quickly whether a given tuple $Z=(z_1,\dots,z_s)$ is suitable.

In other words, given an ideal $I$ in a polynomial ring $K[X]=K[x_1,\dots,x_n]$,
we want to test whether $I\cap K[X\setminus Z]$ can be computed using elimination by substitution.
For this, we want to check whether polynomials of the form $f_i=z_i - h_i$ exist in~$I$,
where $h_i$ is not divisible by any indeterminate in~$Z$. In full generality, this would
require the calculation of a Gr\"obner basis of~$I$ which is infeasible in many cases.
Therefore we search for the polynomials $f_i$ in suitable subspaces of~$I$.
 
The algorithms we develop are very fast and allow us to treat examples with dozens,
even hundreds of indeterminates. The main idea is not to look for the polynomials~$f_i$
directly, but for a tuple of weights~$W$ such that, for any term ordering~$\sigma$ compatible
with the grading given by~$W$, there exist polynomials $f_i\in I$ with $\LT_\sigma(f_i)=z_i$.
Having found~$W$, the actual polynomials defining the $Z$-separating re-embedding can then be
calculated by simple interreduction.

Applications include examples of ideals for which traditional methods of elimination computation
fail, such as ideals defining border basis schemes in Algebraic Geometry. Another area which may
benefit from the new techniques is Cryptography. To apply the new method in this case, we 
developed a variant which is able to perform elimination by substitution for Boolean polynomials.
For example, we show how $Z$-separating Boolean polynomials
speed up several types of algebraic attacks on the cipher AES-128.

%
%

\section*{Acknowledgements}
The first author gratefully acknowledges \textit{Cusanuswerk e.V.}~for financial support.
The third author is supported by the \textit{Vietnam Ministry of Education and Training} 
under the grant number B2025-DHH-02 
and he thanks the University of Passau for its hospitality
during part of the preparation of this paper.

%
%

\end{document}